\documentclass{article}
\usepackage[a4paper, total={6.5in, 8.5in}]{geometry}
\usepackage{graphicx}
\usepackage{amsmath}
\usepackage{amsthm}
\usepackage{xcolor}
\usepackage{float}

\usepackage{algorithm}
\usepackage{algpseudocode}
\usepackage{url}
\usepackage{multirow}
\usepackage{longtable}
\usepackage{setspace}
\usepackage{authblk}

\title{Pickup \& Delivery with Time Windows and Transfers: combining decomposition with metaheuristics}
\author[1]{Ioannis Avgerinos}
\author[1]{Ioannis Mourtos}
\author[1]{Nikolaos Tsompanidis}
\author[1, 2]{Georgios Zois}
\affil[1]{\small ELTRUN Lab, Department of Management Science and Technology, Athens University of Economics and Business,\protect\\ \small Patision 76, 10434, Athens, Greece E-mail: {\tt iavgerinos@aueb.gr, mourtos@aueb.gr, tsompanidisnik@gmail.com, georzois@aueb.gr}}
\affil[2]{\small {Optiscale, Athens 114 72, Greece}}
\date{}

\begin{document}

\maketitle

\begin{abstract}
    This paper examines the generalisation of the Pickup and Delivery Problem that allows mid-route load exchanges among vehicles and obeys strict time-windows at all locations. We propose a novel Logic-Based Benders Decomposition (LBBD) that improves optimality gaps for all benchmarks in the literature and scales up to handle larger ones. To tackle even larger instances, we introduce a refined Large Neighborhood Search (LNS) algorithm that improves the adaptability of LNS beyond case-specific configurations appearing in related literature. 
    
    To bridge the gap in benchmark availability, we develop an instance generator that allows for extensive experimentation. For moderate datasets (25 and 50 requests), we evaluate the performance of both LBBD and LNS, the former being able to close the gap and the latter capable of providing near-optimal solutions. For larger instances (75 and 100 requests), we recreate indicative state-of-the-art metaheuristics to highlight the improvements introduced by our LNS refinements, while establishing its scalability.
\end{abstract}
\textbf{Keywords.} pickup and delivery problem, integer programming, logic-based benders decomposition, large neighbourhood search. 

\section{Background} \label{sec:intro}

\noindent\textbf{Motivation.} The widespread applicability of routing problems has driven the development of numerous variants, along with methods capable of handling large-scale instances. We investigate here a variant that enjoys increasing applicability yet remains quite challenging. Our aim is to provide a novel modelling and decomposition approach, coupled with effective exact and heuristic methods.
    
The traditional \emph{Vehicle Routing Problem} (VRP) \cite{Toth02} assumes vehicles departing from a single depot to cover a single type of operation for all customers (i.e., requests), i.e., either pickups or deliveries. To align with industry requirements, extensions of VRP also consider multiple depots, time windows and combined pickup and delivery operations. Specifically, the \emph{Dial-A-Ride Problem} (DARP) \cite{DARP} assumes pickup and delivery locations per request and hence the precedence constraint that a vehicle must visit the pickup location of a request before its delivery location. In its original version, DARP assumes unitary loads for all requests, as its motivation comes from passenger transport \cite{Sav95}. Therefore, the \emph{Pickup and Delivery Problem} (PDP) \cite{Dum91} extends DARP to accommodate loads of varying size, hence forming the most generic VRP variant. 

It becomes simple to see that VRP reduces to PDP by introducing one copy of the depot per request to serve as its pickup point. Also, PDP allows not just for a different depot per vehicle (as in the multi-depot variant of VRP), but for a different origin and destination depot per vehicle. What further complicates matters is that a vehicle's load fluctuates throughout the route as pickups and deliveries intertwine. The above, combined with the precedence constraints, result in PDP being tackled for instances of minuscule size compared to the scale claimed currently for VRP \cite{Acc24}. 

Further challenges show up in urban logistics and particularly in areas having frequent and densely-concentrated requests. Routing must simultaneously handle strict time windows, restricted access for larger vehicles and rising inventory costs hence lack of warehousing facilities. This explains the appearance of innovative approaches like mobile depots \cite{Hof21} serving as intermediate points for load replenishment, i.e., for transfers of requests between vehicles, or crowdshipping \cite{Voi22} where multiple drivers coordinate to fulfill a request. Therefore, PDP naturally generalises to the \emph{PDP with Transfers} \cite{Cort10}  by introducing locations where loads can be exchanged between vehicles. Such transfer points impose the so called \emph{vehicle synchronisation} constraints: a vehicle picking up a request must arrive at the transfer point after the vehicle dropping it off. Adding time windows for \emph{all} locations (origin and destination of each vehicle, pick and delivery of each request, 
transfer points) yields the \emph{Pickup and Delivery Problem with Time Windows and Transfers}, called hereafter PDPT for conciseness of our notation.

PDPT generalises the Vehicle Routing Problem with Transfers (VRPT) \cite{Drexl13}. Although VRPT incorporates synchronisation constraints, the existing literature considers uniform demand per location rather than distinct pickup and delivery requests. This assumption enables load splitting and simplifies modelling through flow constraints, hence yielding a simpler problem. VRPT has been addressed through Column Generation \cite{Bald17}, and a two-index MILP formulation combined with a Multi-Start Local Search \cite{Agu25}. The study of \cite{Zhang22} incorporates several key aspects of the PDPT, notably multiple depots that can function as transfer points for the deployed vehicles. 

The PDPT imposes several big-M constraints to a monolithic MILP, e.g., precedence constraints ensuring that each request is picked up before being delivered or time-window constraints regulating vehicle arrival times at each location. Also, the inclusion of transfer points inflates the solution space and makes it non-trivial to enforce 
synchronisation constraints. Furthermore, capacity constraints must be managed at every route segment, since pickups and deliveries occur throughout a vehicle's journey. Enhancements like valid inequalities strengthening the lower bound or the elimination of some binary variables remain insufficient even for moderate instances. Current exact methods struggle to produce feasible solutions for up to 30 requests \cite{Lyu23} even without time windows.

Therefore, we have developed a novel Logic-Based Benders Decomposition (LBBD) that reduces optimality gaps for the instances tackled by exact methods \cite{Lyu23} and, most importantly, allows for larger and more complex instances to be tackled. For instance sizes where the efficiency of our exact method diminishes, we introduce a Large Neighbourhood Search (LNS) method that avoids instance-specific characteristics and parameter tuning. Using the LNS solution as a warm-start allows LBBD to offer lower bounds and often improve that solution. \\ 

\noindent\textbf{Literature review. }Mathematical models for the PDP appear in \cite{Dum91, Sav95}, while comprehensive reviews on the evolution of these models appear in \cite{DARP} and in \cite{Koc20}. Due to the complexity of the synchronisation constraints, monolithic MILPs handle instances with very few requests, hence motivating the use of decomposition approaches. An early attempt to accelerate convergence to optimality through decomposition is the Branch-and-Cut scheme of \cite{Cort10}. Although the decomposition outperforms the monolithic MILP, instances with only up to 6 requests, in a specific DARP scenario involving passenger transport (i.e., with unitary loads) have been solved. 
The improved MILP formulation of \cite{Rais14} expands the solvable instance size to 14 requests. The approach of \cite{Mass14}, combining Column Generation and Branch-and-Cut, performs well on instances with up to 193 requests but for the special PDPT case where routes of predefined links (a.k.a. shuttle routes) serve passengers for a limited number of stops.
    
The most recent advancement in mathematical modelling for PDPT is presented by \cite{Lyu23} that extends \cite{Rais14}. This new MILP addresses instances of up to 30 requests with loose time windows but only up to 5 requests when time windows become strict. Given the absence of publicly available benchmark instances of larger size, we rely upon the work of \cite{Lyu23} to assess our methods. 

Several studies have explored the use of transshipment (i.e., transfers) in PDP and its variants by deploying various metaheuristics. Early work \cite{Mitro} investigates the benefits of load transfers between vehicles in courier operations and particularly in large, multi-city areas. Building on this, two nearly-concurrent studies pioneer the development of Adaptive Large Neighbourhood Search (ALNS) for the PDPT: \cite{Qu} proposes a GRASP approach that employs ALNS to improve the initially constructed solutions, while \cite{Masson2013-ae} proposes a handful of efficient heuristics for modifying a partial solution and then incorporating it into an ALNS scheme. \cite{Danloup2018-lh} are the first to examine the potential of Genetic Algorithms (GA), hence proposing both an LNS and a GA. Finally, transfers have been examined in crowd-sourced urban freight delivery \cite{Sam20}, where ALNS effectively optimises transfer opportunities and vehicle synchronisation, while significantly reducing travel distances and the number of required drivers. These findings have motivated us to develop a refined LNS scheme, whose core aim is to complement the strengths of our LBBD approach.\\

\noindent \textbf{Scope and contribution}. We provide an effective combination of exact and metaheuristic methods for PDPT, hence consolidating a currently fragmented literature. Our LBBD is a more scalable exact method, as prior work \cite{Cort10, Rais14, Lyu23} shows that synchronisation constraints, together with the large number of variables, severely limits the effectiveness of monolithic MILPs. We employ a strong relaxation serving as a master problem, which constructs parallel paths from the pickup to the delivery point of all requests. Most properties of the PDPT - the capacity constraints and precedences of pickup and delivery operations - are satisfied with significantly fewer variables, thus leaving vehicle assignment to routes plus vehicle synchronisation to an easy to solve subproblem. The use of additional valid inequalities and standard combinatorial (i.e., Benders') cuts improve the convergence to optimality.


Our LNS-based method emphasises simplicity, robustness, and minimal parameter tuning. Building on the LNS framework of \cite{Pisinger2019-wh}, our key contributions include a parameter-free destroy operator 
for adaptive request removal, inspired by \cite{Danloup2018-lh}, a normalised insertion difficulty metric \cite{Danloup2018-lh} that enhances the repair process without instance-specific tuning, and a combination of a randomised cheapest-insertion heuristic \cite{Christiaens2020-vb} with constant-time feasibility checks \cite{Masson2013-xm}. 
Together with a delayed acceptance mechanism \cite{Burke2017-vf} that improves solution quality, these yield an effective approach for PDPT, either to support LBBD or as standalone for large instances.

Given the scarcity of PDPT benchmarks, we adapt datasets from more restricted variants. Let us recall that benchmarks from \cite{Lyu23} have up to 30 requests, 3 transfer points, and 3 vehicles for loose time windows and only up to 5 requests for strict time windows. 
The most extensive dataset from the metaheuristics' community covers up to 100 requests \cite{Sam20}, yet these instances are too specific for crowd-sourced drivers (e.g., they assume unlimited number of vehicles and neutralised capacity constraints).
    
To bridge the gap in benchmark availability, we develop a new instance generator inspired by \cite{Sar20}. Our generator introduces more robust scenarios by varying the tightness of time windows and incorporating transfer point construction. The generated instances present significant challenges for an exact method, whereas the LNS algorithm consistently produces high-quality solutions. To evaluate these solutions, we use the LBBD scheme as a lower bound provider, demonstrating that the LNS solutions accomplish an up to 40\% optimality gap. The LBBD struggles to find feasible solutions when used alone (under a 1-hour time limit), but speeds up once starting from the LNS solution (which it often improves). \\

       
\noindent \textbf{Outline.} This paper describes a comprehensive toolbox of exact and metaheuristic methods for the PDPT, hence its remainder is organised as follows. Section \ref{sec:prel} overviews our methodological framework and formally introduces PDPT along with our notation. Section \ref{sec:exact} outlines the intuition of our LBBD and details the decomposition scheme. Section \ref{sec:lns} elaborates on our refined LNS. The benchmark generator and experimental results are discussed in Section \ref{sec:experiments}. Section \ref{sec:conclusions} summarises our findings to motivate further work.
\section{Preliminaries} \label{sec:prel}

\subsection{Logic-Based Benders Decomposition}
LBBD \cite{hooker03} is an extension of the classical Benders Decomposition \cite{benders}. Let $\mathcal{P}$ be a minimisation problem of two groups of variables $x$, $y$ in domains $D_{x}$ and $D_{y}$: 
        \begin{equation}
            \mathcal{P}:\qquad min\{f(x) + g(y):C_{1}(x), C_{2}(y), C_{3}(x, y), x\in D_{x}, y\in D_{y}\}. \notag 
        \end{equation}
Cost functions $f(x)$ and $g(y)$ are assumed to be non-negative. The problem is subject to constraints $C_{1}(x)$, $C_{2}(y)$ and $C_{3}(x, y)$. 
To partition $\mathcal{P},$ $\mathcal{M}$ is defined as
        \begin{equation}
            \mathcal{M}:\qquad min\{z:C_{1}(x), x\in D_{x}, z\geq f(x)\}. \notag 
        \end{equation}
        
\noindent As $f(x)$ and $g(y)$ are non-negative, $z$ is a lower bound of the objective function of $\mathcal{P}.$ Also all constraints of $\mathcal{M}$ are also imposed to $\mathcal{P}$. Therefore, $\mathcal{M}$ is a \emph{relaxation} of $\mathcal{P}$ called the \emph{master problem}.

Given the optimal solution $\bar{x}$ of $\mathcal{M}$ and the respective objective value $\bar{z}$, the following problem $\mathcal{S}$ computes a local optimum of $\mathcal{P}$:
        \begin{equation}
            \mathcal{S}:\qquad min\{f(\bar{x}) + g(y):C_{2}(y), C_{3}(\bar{x}, y), y\in D_{y}\}. \notag 
        \end{equation}
We notice that $\mathcal{S}$ is a \emph{subproblem} of $\mathcal{P}$ that finds the optimal solution in the restricted search space defined by $\bar{x}$. If $C_{3}(\bar{x}, y)$ are satisfied for the given solution $\bar{x}$, the optimal objective value of $\mathcal{S},$ denoted as $\bar{\zeta},$ is an upper bound of the optimal objective value of $\mathcal{P}$. Hence, $\frac{\bar{\zeta} - \bar{z}}{\bar{\zeta}}$ defines a valid optimality gap.

That is, LBBD employs an iterative exchange of knowledge between $\mathcal{M}$ and $\mathcal{S}$: $\mathcal{M}$ provides a lower bound of $\mathcal{P}$ and a partial solution $\bar{x}$, and, then, $\mathcal{S}$ returns an upper bound which is the local optimum of the restricted search space. A valid combinatorial cut ensures that the master problem avoids repeating the same partial solution at any future iteration. Eventually, the lower bound and the best found upper bound converge to identical values, thus reaching the global optimal solution of $\mathcal{P}$.

If $\mathcal{S}$ can provide a feasible solution, then the combinatorial cut is a logic constraint of the following form, called \emph{optimality cut}:
        \begin{equation}
            \text{if }x = \bar{x} \rightarrow z\geq \bar{\zeta}. \notag
        \end{equation}
        On the contrary, if $\mathcal{S}$ is infeasible for a partial solution $\bar{x}$, then a logic constraint called \emph{feasibility cut} is generated:
        \begin{equation}
            x \neq \bar{x}.\notag
        \end{equation}
Both types of logic constraints can be linearised trivially, if variables $x$ are binary.

A hybridisation of Branch-and-Cut algorithms with LBBD, namely \emph{Branch-and-Check} \cite{thor}, has shown faster convergence, specially for partitions that involve small subproblems. In this case, the master problem is not solved optimally in an iterative manner; instead, the branching tree is explored once, interrupting the solution at each integer node to retrieve the local upper bound or prune the branch in case of infeasibility.
            
\subsection{Large Neighbourhood Search} 
LNS \cite{Pisinger2019-wh} is an iterative local search method that explores a large neighborhood of a given solution in each iteration. This neighborhood is implicitly defined through two solution modification phases: first, the solution is partially destroyed by removing a subset of the requests and, then, the removed requests are reinserted to construct a new solution. The new solution is evaluated against the current one, and accepted based on a predefined criterion.
        
Hence, the main LNS components are the destroy operator defining how requests are removed from given solution, the repair operator determining how the removed requests are reinserted, and the acceptance criterion governing whether the repaired solution is accepted. A stopping criterion is also necessary as LNS cannot guarantee convergence to optimality. The design of these components influences the exploration of the solution space, impacting both the efficiency and effectiveness of the method.
        
\begin{algorithm}
            \caption{The LNS metaheuristic}
            \label{alg:lns}
            \begin{algorithmic}[1]
            \Function{LNS}{$G$}
                \State $s \leftarrow$ \Call{Initial-Solution}{$G$}
                \State $s^* \leftarrow s$
                \While{\textit{stopping criterion not met}}
                    \State $s' \leftarrow$ \Call{Local-Search}{$s$}
                    \If{\Call{Accept}{$s'$}}
                        \State $s \leftarrow s'$
                    \EndIf
                    \If{$cost(s') < cost(s^*)$}
                        \State $s^* \leftarrow s'$
                    \EndIf
                \EndWhile
                \State \textbf{return} $s^*$
            \EndFunction
            \end{algorithmic}
            \end{algorithm}

Algorithm \ref{alg:lns} presents the LNS workflow in an abstract manner. The method begins with an initial solution $s$ for the given problem instance, which is also stored as the best-known solution $s^*$. The algorithm then enters the local search process, which continues until a termination criterion is met. In each iteration, the neighborhood of the current solution $s$ is explored by sequentially applying a pair of destroy and repair operators. The resulting solution $s'$ replaces the current solution, if it meets the acceptance criterion. Additionally, if $s'$ has a lower cost than $s^*$, it is stored as the best solution found so far.

\subsection{Notation}
The PDPT is defined on a set of requests $\mathcal{R}$ for pick-up and delivery, a set of locations $\mathcal{J}$ and a fleet of vehicles $\mathcal{K}$. The set of locations $\mathcal{J}$ comprises all depots $\mathcal{D}$, transfer points $\mathcal{T}$ and demand points $\mathcal{C}$ (both pick-up and delivery), the three subsets not being necessarily disjoint. Each request $r\in \mathcal{R}$ has a demand $q_{r}$ (e.g., weight, volume) and must be picked up from a specific location $p_{r} \in \mathcal{C}$ and delivered to another location $d_{r} \in \mathcal{C}$. The requests are fulfilled by a fleet of homogeneous vehicles $\mathcal{K}$, each with capacity $Q$. Each vehicle $k\in \mathcal{K}$ begins its route at an origin depot $o_{k} \in \mathcal{D}$ and completes it at a destination depot $e_{k} \in \mathcal{D}$.
    
        \begin{table}[tbh]
            \scriptsize
            \centering
            \caption{Sets and parameters}
            \label{tab:annotations}
            \begin{tabular}{lll} \\
            \hline
            \multicolumn{3}{l}{\textbf{Sets}}                                                        \\ \hline
            $\mathcal{R}$          & \multicolumn{2}{l}{Requests}                                          \\
            $\mathcal{J}$          & \multicolumn{2}{l}{Locations}                                         \\
            $\mathcal{K}$          & \multicolumn{2}{l}{Vehicles}                                          \\
            $\mathcal{T}$          & \multicolumn{2}{l}{Transfer points; $\mathcal{T}\subset \mathcal{J}$} \\
            $\mathcal{D}$          & \multicolumn{2}{l}{Depots; $\mathcal{D}\subset \mathcal{J}$} \\
            $\mathcal{C}$          & \multicolumn{2}{l}{Demand points; $\mathcal{C}\subset \mathcal{J}$} \\\hline
            \multicolumn{3}{l}{\textbf{Parameters}}                                                      \\ \hline
            $p_{r}\in \mathcal{J}$ & $r\in \mathcal{R}$         & Pickup point of request $r$              \\
            $d_{r}\in \mathcal{J}$ & $r\in \mathcal{R}$         & Delivery point of request $r$            \\
            $q_{r}$                & $r\in \mathcal{R}$         & Demand of request $r$                      \\
            $Q$                    &                            & Capacity of vehicles                     \\
            $o_{k}\in \mathcal{J}$ & $k\in \mathcal{K}$         & Origin depot of vehicle $k$              \\
            $e_{k}\in \mathcal{J}$ & $k\in \mathcal{K}$         & Destination depot of vehicle $k$         \\
            $l_{j}$                & $j\in \mathcal{J}$         & Earliest time of visit at $j$            \\
            $u_{j}$                & $j\in \mathcal{J}$         & Latest time of visit at $j$  \\
            $t_{ij}$               & $i, j\in \mathcal{J}$      & Travel time from $i$ to $j$ \\
            $c_{ij}$               & $i, j\in \mathcal{J}$      & Distance from $i$ to $j$ \\\hline
            \end{tabular}
        \end{table}
Notably, when a specific request $r$ is loaded onto a vehicle $k$, it can then be unloaded at a subsequent stop by the same or a different vehicle. In the former case, both $p_{r}$ and $d_{r}$ are included in the same route, ensuring that $p_{r}$ precedes $d_{r}$. In the second case, vehicles may utilise designated a transfer point in $\mathcal{T}$, where requests can be offloaded and, subsequently, picked up by other vehicles for delivery or further transport.
    
Each location $j\in \mathcal{J}$ is associated with a specific time window $[l_{j}, u_{j}]$, which defines the permissible time interval during which the location can be visited for pickup or delivery. Each pair of locations $(i, j)$ has a travel time $t_{ij}$ and a distance $c_{ij}$. For simplicity we assume zero service times, i.e., pickups and deliveries are instantaneous. This is at no loss of generality, as service time depend on the location alone hence can be added to the travel time from any other location.

The objective here is to determine a set of vehicle routes minimising total distance travelled, that serve all requests while adhering to the specified time windows and vehicle capacity constraints. Each route must start at some vehicle's origin depot and end at its designated destination depot. 

\section{Benders decomposition} \label{sec:exact}
    \subsection{An improved exact approach}

Let us recall from \cite{Cort10}, for self-sufficiency of our presentation, the following five key properties of PDPT: (pdpt1) each route must begin at the vehicle's origin depot and end at its destination depot without forming subtours, (pdpt2) each pickup and delivery point must be visited exactly once, (pdpt3) pickups must always precede their corresponding deliveries, (pdpt4) load transfers between vehicles must be synchronised (i.e., the vehicle unloading a request must arrive at the transfer location before the vehicle loading it), and (pdpt5) capacity constraints must be maintained throughout a vehicle's route. As we move beyond \cite{Cort10} we must introduce property (pdpt6): all visits must comply with predefined time windows. We refer to these properties in the sequel to explain the correctness of our mathematical modelling.

PDPT has been formulated as a MILP in \cite{Rais14} and \cite{Lyu23}. 
A key limitation lies in the MILP's reliance on four-indexed binary variables $y_{rijk}$, which are set to 1 if request $r$ is loaded onto vehicle $k$ while traveling from location $i$ to $j$. Given that the total number of locations $\mathcal{J}$ includes pickup and delivery points for each request, origin and destination depots for each vehicle, and all transfer points ($2\cdot |\mathcal{R}|+2\cdot |\mathcal{K}|+|\mathcal{T}|$ locations), the number of such variables explodes even for moderate instances.

Partitioning the MILP can facilitate the computation of feasible solutions. The work of \cite{Cort10} already considers a relaxation of the holistic model by omitting complicating constraints, and a feasibility subproblem to validate if the obtained partial solution can lead to a feasible routing schedule. Combinatorial cuts \cite{Codato06} are added to the relaxation, if a feasible solution is not obtained. This decomposition has not shown increased scalability, but accelerates the optimal solution of small-scaled instances (up to 6 requests). 

We therefore introduce a novel partitioning strategy, leading to a stronger relaxation of the original problem as the master problem. Specifically, we exclude from the master problem the set of vehicles $\mathcal{K}$ and the synchronisation of transfers (pdpt4). The intuition behind this choice is to enforce properties (pdpt2), (pdpt3), and (pdpt5) while relaxing (pdpt1) and (pdpt6) in a way that increases the likelihood of generating a feasible routing plan. The solution to the master problem consists of a set of parallel paths, each corresponding to a request $r\in \mathcal{R}$. These paths inherently satisfy (pdpt2), (pdpt3), and (pdpt5), thus allowing a simpler subproblem to focus on assigning selected path edges to valid routes, while addressing the remaining properties.

\begin{figure}[H]
           \centering \includegraphics[scale=0.4]{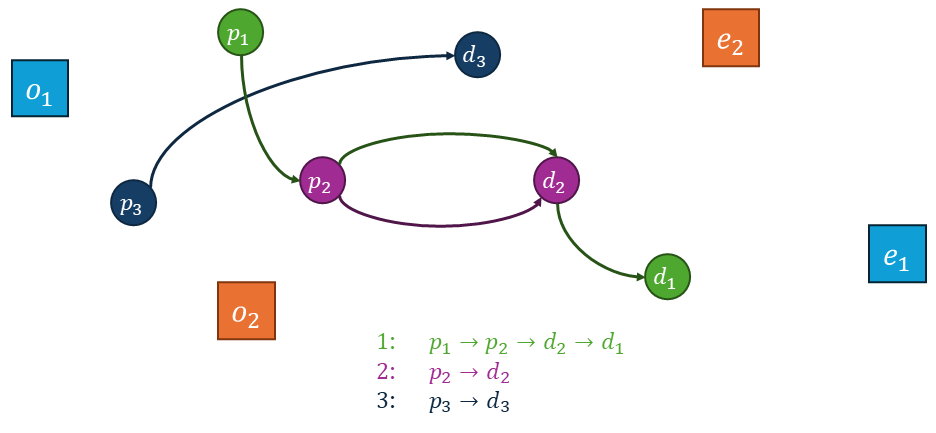}
           \caption{Example of path construction for three requests and two vehicles}\label{fig:subpaths}
\end{figure}

Figure \ref{fig:subpaths} illustrates an example involving three requests of size 1 and two vehicles with a capacity of 2. Each request $r$ has each own path, beginning from $p_{r}$ and ending at $d_{r}$. The paths of requests 1 and 2 overlap, as they share trip $p_{2}\rightarrow d_{2}$ (hence the two parallel edges), whereas the path for request 3 is disjoint. Importantly, incorporating this partial solution into the original problem, while preserving the path edges, significantly simplifies the optimisation problem. Additionally, several constraints of the original problem are inherently satisfied by this partial solution: vehicle capacities are not exceeded for any path, and the precedence of pickup and delivery operations is maintained. Such a construction of paths is clearly a relaxation of the original problem, as it omits several constraints, most notably the involvement of vehicles. 


\subsection{The master problem}

\textbf{The master MILP.} The formulation relies on two sets of binary variables. Variables $y_{rij}$ are set to $1$ if request $r\in \mathcal{R}$ is loaded (on some vehicle) during trip (i.e., edge) $(i, j),$ where $i,j\in \mathcal{J}$, and $0$ otherwise. Variables $x_{ij}$ indicate whether trip $(i, j)$ is selected or not. Let us now explain why the above formulation is correct by associating constraints with the properties (pdpt1)-(pdpt6) listed at the end of Section 2.1.
\begin{flalign}
            \text{min } & z && \notag &&\\
                        & z\geq \sum_{i\in \mathcal{J}}\sum_{j\in \mathcal{J}}c_{ij}\cdot x_{ij} && \label{eq:m1} &&\\
                        &\sum_{j\in \mathcal{J}}y_{rp_{r}j} = 1 && r\in \mathcal{R} \label{eq:m2} &&\\
                        &\sum_{j\in \mathcal{J}}y_{rjd_{r}} = 1 &&  r\in \mathcal{R} \label{eq:m3} &&\\
                        &\sum_{i\in \mathcal{J}}y_{rij}\leq 1 &&  r\in \mathcal{R}, j\in \mathcal{J}\notin \{p_{r}, d_{r}\} \label{eq:m4} &&\\
                        &\sum_{i\in \mathcal{J}}y_{rij} = \sum_{i\in \mathcal{J}}y_{rji} && r\in \mathcal{R}, j\in \mathcal{J}\notin \{p_{r}, d_{r}\} \label{eq:m5} &&\\
                        &\sum_{i\in \mathcal{J}}x_{ij} = 1 && j\in \mathcal{C} \label{eq:m6} &&\\
                        &\sum_{i\in \mathcal{J}}x_{ij} = \sum_{i\in \mathcal{J}}x_{ji} && j\in \mathcal{J}\setminus \mathcal{D} \label{eq:m7} &&\\
                        &\sum_{r\in \mathcal{R}}q_{r}\cdot y_{rij} \leq Q\cdot x_{ij} && i\in \mathcal{J}, j\in \mathcal{J} \label{eq:m8} &&\\
                        &\text{[Valid inequalities]} &&\label{eq:m9} &&\\
                        &\text{[Time windows constraints]} &&\label{eq:m10} &&\\
                        &\text{[Benders cuts]} &&\label{eq:m11} &&\\
                        &x_{ij}\in \{0, 1\} && i\in \mathcal{J}, j\in \mathcal{J} \notag &&\\
                        &y_{rij}\in \{0, 1\} && r\in \mathcal{R}, i\in \mathcal{J}, j\in\mathcal{J}, i\neq j \notag &&     
        \end{flalign}
        
The objective function is the total distance travelled (\ref{eq:m1}). Constraints (\ref{eq:m2}) and (\ref{eq:m3}) ensure that each request is loaded on some edge departing from its pickup point and on some edge arriving at its delivery point; notice that this also ensures that property (pdpt3) is eventually satisfied. Constraints (\ref{eq:m4}) guarantee that each request visits any location at most once; if a transfer point is eligible for multiple visits, the set $\mathcal{T}$ should include duplicate copies of that location, as also noted in \cite{Cort10}. Constraints (\ref{eq:m5}) preserve the flow, i.e., they ensure that each request departs from any visited location, except for its pickup and delivery points.
        
Regarding now variables $x_{ij}$, Constraints (\ref{eq:m6}) enforce that each demand point is visited once, while Constraints (\ref{eq:m7}) ensure flow conservation at non-depot locations: these constraints satisfy property (pdpt2), while also enforcing (pdpt1). Constraints (\ref{eq:m8}) serve a dual purpose: they enforce that each trip carries a request volume within the vehicle capacity and establish the connection between variables $y_{rij}$ and $x_{ij}$. Since the vehicles have identical capacity, Constraints (\ref{eq:m8}) satisfy property (pdpt5). \\

\noindent \textbf{Valid inequalities. }Constraints (\ref{eq:m2})–(\ref{eq:m8}) capture most of the PDP intricacies. Notice, however, that Constraints (\ref{eq:m6}) do not involve transfer points and allow multiple incoming and outgoing flows through locations in $\mathcal{T}$. To ensure the validity of paths from $p_{r}$ to $d_{r}$ for all requests $r\in \mathcal{R}$, eliminate invalid trips, and strengthen the MILP lower bound, we also incorporate the following set of valid inequalities.
\begin{flalign}
            &x_{jj} = 0 && j\in \mathcal{J} \label{eq:v1} &&\\
            &x_{d_{r}p_{r}} = 0 && r\in \mathcal{R} \label{eq:v2} &&\\
            &x_{o_{k}d_{r}} = 0 && k\in \mathcal{K}, r\in \mathcal{R} \label{eq:v3} &&\\
            &x_{p_{r}e_{k}} = 0 && k\in \mathcal{K}, r\in \mathcal{R} \label{eq:v4} &&\\
            &x_{o_{k}e_{k'}} = 0 && k\in \mathcal{K}, k'\mathcal{K}, k\neq k' \label{eq:v5} &&\\
            &x_{jo_{k}} = 0 && k\in \mathcal{K}, j\in \mathcal{J} \label{eq:v6} &&\\
            &x_{e_{k}j} = 0 && k\in \mathcal{K}, j\in \mathcal{J} \label{eq:v7} &&\\
            &\sum_{i\in \mathcal{J}}y_{rij} = 0 && r\in \mathcal{R}, j\in \mathcal{D} \label{eq:v8}&&\\
            &\sum_{j\in \mathcal{J}}y_{rd_{r}j} = 0 && r\in \mathcal{R} \label{eq:v9}&&\\
            &\sum_{j\in \mathcal{J}}y_{rjp_{r}} = 0 && r\in \mathcal{R} \label{eq:v10}&&\\
            &\sum_{j\in \mathcal{J}}x_{o_{k}j} = 1 && k\in \mathcal{K} \label{eq:v11} &&\\
            &\sum_{j\in \mathcal{J}}x_{je_{k}} = 1 && k\in \mathcal{K} \label{eq:v12} &&
        \end{flalign}

Constraints (\ref{eq:v1})–(\ref{eq:v10}) eliminate a significant number of variables, thereby reducing the search space. Specifically, Constraints (\ref{eq:v1}) prohibit loops, while (\ref{eq:v2}) prevent trips from a delivery point back to the corresponding pickup point. Constraints (\ref{eq:v3}) and (\ref{eq:v4}) ensure that no delivery point is the first stop, and no pickup point is the last stop of any vehicle.

Although direct connections between an origin and a destination depot are not explicitly forbidden - since such trips indicate that a vehicle remains unused - the connection of an origin depot with the destination depot of a different vehicle is invalid, hence imposing (\ref{eq:v5}). Origin depots cannot serve as destinations, as imposed by (\ref{eq:v6}), and destination depots cannot serve as origins, as dictated by (\ref{eq:v7}). Depot locations must not be included in the path of any request, as stated by (\ref{eq:v8}).
        
Constraints (\ref{eq:v9}) and (\ref{eq:v10}) ensure, respectively that requests cannot be loaded before reaching their pickup point and they remain loaded until reaching their delivery point. Furthermore, (\ref{eq:v11}) and (\ref{eq:v12}) guarantee that depots are incorporated into the routing plan.

The following inequalities reduce the likelihood of partial solutions that lead to infeasibility:
        \begin{flalign}
            &\sum_{i\in \mathcal{J}}x_{ij}\leq |\mathcal{K}| && j\in \mathcal{T}. \label{eq:v13} &&
        \end{flalign}
\noindent Let us explain this. As (\ref{eq:m6}) do not apply to $\mathcal{T}$, there are no restrictions on assigning transfer points to trips. However, each vehicle can visit a transfer point at most once, thus we should prevent an overuse of edges accessing the same transfer point. Constraints (\ref{eq:v13}) limit the number of visits per transfer point to $|\mathcal{K}|$, hence ensuring that the selected edges can be validly assigned to the vehicles in the subproblem presented next.
        
We note that limiting the number of visits to one per (vehicle, transfer point) pair is a modeling decision rather than a strict constraint, as transfer points might be used without such a restriction. To accommodate multiple visits, the set $\mathcal{T}$ should include duplicate copies of transfer points, allowing them to be visited more often (a similar approach is suggested by \cite{Cort10}).\\

\noindent \textbf{Time windows constraints. } 
Time windows impose additional constraints that account for travel times along selected trips. Therefore, we introduce two sets of continuous variables. The variables $a_{j}$ represent the arrival times at all locations except transfer points. Since transfer points can be visited by multiple vehicles, their arrival times cannot be uniquely defined. Therefore, we introduce variables $b_{rj}$, which indicate the arrival time of request $r$ at transfer point $j$. All arrival time variables are constrained within the time window limits $l_{j}$ and $u_{j}$. The constraints are as follows.
\begin{flalign}
            &a_{i} + t_{ij} - a_{j}\leq (t_{ij} + u_{i})\cdot (1 - x_{ij}) && \forall i\in \mathcal{J}\setminus \mathcal{T}, j\in \mathcal{J}\setminus \mathcal{T} \label{eq:t1} &&\\
            &a_{i} + t_{ij} - b_{rj}\leq (t_{ij} + u_{i})\cdot (1 - y_{rij}) && \forall r\in \mathcal{R}, i\in \mathcal{J}\setminus \mathcal{T}, j\in \mathcal{T} \label{eq:t2} &&\\
            &b_{ri} + t_{ij} - a_{j}\leq (t_{ij} + u_{i})\cdot (1 - y_{rij}) && \forall r\in \mathcal{R}, i\in \mathcal{T}, j\in \mathcal{J}\setminus \mathcal{T} \label{eq:t3} &&\\
            &b_{ri} + t_{ij} - b_{rj}\leq (t_{ij} + u_{i})\cdot (1 - y_{rij}) && \forall r\in \mathcal{R}, i\in \mathcal{T}, j\in \mathcal{T} \label{eq:t4} &&\\
            &l_{j}\leq a_{j}\leq u_{j} && j\in \mathcal{J} \notag &&\\
            &l_{j}\leq b_{rj}\leq u_{j} && r\in \mathcal{R}, j\in \mathcal{T} \notag &&
        \end{flalign}
\noindent Constraints (\ref{eq:t1}) determine the arrival times of non-transfer points (interestingly, similar constraints are commonly used for computing completion times of operations which are determined by sequence-dependent times). Constraints (\ref{eq:t2})-(\ref{eq:t4}) define the arrival times of trips which involve transfer points: since each transfer point can be visited for more than once, a single-index variable $a_{t}$ for each $t\in \mathcal{T}$ would be invalid.
\subsection{The subproblem}
        \textbf{Preprocessing. }For any solution of the master problem, let $\bar{y}_{rij}$ be the values of variables $y_{rij}$. The following simple procedure constructs the path of each request $r\in \mathcal{R}$:

        \begin{algorithm}[H]
        	\fontsize{9}{11}\selectfont
        	\caption{\label{alg:paths}Construction of path for request $r$}
            \begin{algorithmic}[1]
                \State Set $\texttt{path}_{r}$ = \{$p_{r}$\} be the path of $r$ and $\lambda = p_{r}$ be the last visited location of the path;
                \While{$\lambda \neq d_{r}$}
                    \For{$j\in \mathcal{J}$}
                        \If{$\bar{y}_{r\lambda j} = 1$}
                            \State Append $j$ to $\texttt{path}_{r}$;
                            \State $\lambda \leftarrow j$;
                            \State \textbf{break}
                        \EndIf
                    \EndFor
                \EndWhile
                \State Return $\texttt{path}_{r}$;
            \end{algorithmic}
        \end{algorithm}

\noindent Fixing the selected trips in paths leads to a subproblem, in which a valid vehicle routing plan must be generated, while adhering to the paths determined by the master problem. Since capacity constraints and the precedence of pickups and deliveries are already ensured, the focus shifts to synchronising operations at transfer points while respecting time windows for all locations.

To achieve this synchronisation, we construct vectors $\tau_{r}$ for each request $r$, containing triples of consecutive locations in $\texttt{path}_{r}$, where the intermediate location is a transfer point. For example, if trips $(i, t)$ and $(t, j)$ are included in $\texttt{path}_{r}$ and $t$ is a transfer point, then the triple $(i, t, j)$ is appended to $\tau_{r}$. We also define set $\mathcal{E}$, containing all trips which are included in any path: $\mathcal{E} = \{(i, j):\sum_{r\in \mathcal{R}}\bar{y}_{rij}\geq 1\}$.\\

\noindent \textbf{MILP formulation. }The subproblem is formulated as an MILP and aims to assemble the paths from the master problem into complete routes for the available vehicles $\mathcal{K}$. This is done while respecting the time windows of all requests and ensuring synchronisation between drop-offs and pickups at transfer points. The subproblem is required to retain all trips $(i, j)$ included in the paths of requests $\mathcal{R}$ and is permitted to generate additional trips that connect disjoint paths into a single route. 
\begin{flalign}
            \text{min } & \sum_{i\in \mathcal{J}}\sum_{j\in \mathcal{J}}\sum_{k\in \mathcal{K}}c_{ij}\cdot x_{ijk} && \notag &&\\
                        &\sum_{j\in \mathcal{J}}x_{o_{k}jk} = 1 &&  k\in \mathcal{K} \label{eq:s1} &&\\
                        &\sum_{j\in \mathcal{J}}x_{je_{k}k} = 1 && k\in \mathcal{K} \label{eq:s2} &&\\
                        &\sum_{j\in \mathcal{J}}x_{jo_{k}k} = 0 &&  k\in \mathcal{K} \label{eq:s3} &&\\
                        &\sum_{j\in \mathcal{J}}x_{e_{k}jk} = 0 &&  k\in \mathcal{K} \label{eq:s4} &&\\
                        &\sum_{i\in \mathcal{J}}x_{ijk} \leq 1 && j\in \mathcal{J}, k\in \mathcal{K} \label{eq:s5} &&\\
                        &\sum_{i\in \mathcal{J}}x_{ijk} = \sum_{i\in \mathcal{J}}x_{jik} && j\in \mathcal{J}\setminus \mathcal{D}, k\in \mathcal{K} \label{eq:s6} &&\\
                        &\sum_{j\in \mathcal{J}}x_{p_{r}jk} = 1 &&  r\in \mathcal{R}, k\in \mathcal{K} \label{eq:s7} &&\\
                        &\sum_{j\in \mathcal{J}}x_{jd_{r}k} = 1 &&  r\in \mathcal{R}, k\in \mathcal{K} \label{eq:s8} &&\\
                        &\sum_{k\in \mathcal{K}}x_{ijk} = 1 &&  (i, j)\in \mathcal{E} \label{eq:s9} &&\\
                        &a_{ik} + t_{ij} - a_{jk}\leq (t_{ij} + u_{i})\cdot (1 - x_{ijk}) && i\in \mathcal{J}, j\in \mathcal{J}, k\in \mathcal{K} \label{eq:s10} &&\\
                        &a_{p_{r}k}\leq a_{d_{r}k} && r\in \mathcal{R}, k\in \mathcal{K} \label{eq:s11} &&\\
                        &a_{tk} - a_{tl}\leq u_{t}\cdot (2 - x_{itk} - x_{tjl}) && r\in \mathcal{R}, (i, t, j)\in \tau_{r}, k\in \mathcal{K}, l\in \mathcal{K} \label{eq:s12} &&\\
                        &x_{ijk}\in \{0, 1\} && i\in \mathcal{J}, j\in \mathcal{J}, i\neq j, k\in \mathcal{K} \notag &&\\
                        &l_{j}\leq a_{jk}\leq u_{J} && j\in \mathcal{J}, k\in\mathcal{K} \notag &&    
        \end{flalign}

By means of an example, if the solution to the master problem were as illustrated in Figure \ref{fig:subpaths}, the subproblem would be compelled to generate a single route for the paths of Requests $1$ and $2$, as they share the trip $(p_{2}, d_{2})$. Additionally, it would determine whether the path for Request 3 should be incorporated into the same route or assigned to a different vehicle.

Variables $x_{ijk}$ extend the functionality of the relaxed variables $x_{ij}$ from the master problem by incorporating the set of vehicles. If vehicle $k$ includes trip $(i, j)$ in its route, then $x_{ijk}$ is set to 1; otherwise, it remains 0. Additionally, continuous variables $a_{jk}$ represent the arrival times of vehicle $k$ at location $j$.

Each vehicle starts at its origin depot and returns to its destination depot, as enforced by Constraints (\ref{eq:s1})–(\ref{eq:s4}). Each location can be visited at most once, hence imposing (\ref{eq:s5}), while flow conservation is maintained at all non-depot locations as designated by (\ref{eq:s6}). Constraints (\ref{eq:s7}) and (\ref{eq:s8}) ensure that all pickup and delivery points are visited. Moreover, Constraints (\ref{eq:s9}) guarantee that all trips selected in the master problem are included in the solution.
        
Arrival times at visited locations are defined by Constraints (\ref{eq:s10}), while (\ref{eq:s11}) enforce the correct order of pickup and delivery for each request. Finally, (\ref{eq:s12}) ensure vehicle synchronisation at transfer points: a vehicle picking up a request must arrive at the transfer point after the vehicle dropping it off.

\subsection{The Branch-and-Check method}
Recall that $\mathcal{P}$ denotes the original PDPT and $\mathcal{M}$ the master problem in the decomposition. Solving $\mathcal{M}$ generates solutions for iterations $n = 1, 2, \ldots$, each associated with a set of selected edges $\mathcal{E}^{n}$ and a vector $\tau^{n} = \{\tau_{r}^{n}:r\in \mathcal{R}\}$. The corresponding subproblem $\mathcal{S}(\mathcal{E}^{n}, \tau^{n})$ then constructs a complete routing plan, yielding a feasible solution to $\mathcal{P}$ as long as the fixed edges permit it. The objective value of the optimal solution of $\mathcal{M}$ is a lower bound of the objective value of the optimal solution of $\mathcal{P}$, while the objective value of the subproblem $\bar{\zeta^{n}}$ at any iteration $n$ is an upper bound. We solve the master problem $\mathcal{M}$ through Branch-and-Check, i.e., whenever an MILP solver finds an integer solution to $\mathcal{M}$, the solution process is interrupted, and the corresponding subproblem evaluates feasibility. 

If the subproblem is feasible, let $\bar{\zeta}^{n}$ denote its objective value: the following optimality cut ensures that if the master problem attempts to reproduce the same solution, its objective value is constrained to the actual upper bound $\bar{\zeta}^{n}:$
\begin{flalign}
            z\geq \bar{\zeta}^{n} - \bar{\zeta}^{n}\cdot (|\mathcal{E}^{n}| - \sum_{(i, j)\in \mathcal{E}^{n}}x_{ij}). &&  \label{eq:optCut}
        \end{flalign}
\noindent Otherwise, i.e., if the selected edges $\mathcal{E}^{n}$ fail to produce a feasible routing plan, we add to $\mathcal{M}$ the standard feasibility cut
\begin{flalign}
            \sum_{(i, j)\in \mathcal{E}^{n}}x_{ij}\leq |\mathcal{E}^{n}| - 1. &&  \label{eq:feasCut}
        \end{flalign}
After exploring the solution at iteration $n$, the solution process for $\mathcal{M}$ continues. If the latest upper bound improves upon all previously found bounds, the best lower and upper bounds are combined to update the optimality gap. Once the gap reaches zero, the routing plan corresponding to the best upper bound is retrieved, representing the global optimal solution of $\mathcal{P}$.

The validity of (\ref{eq:optCut}) and (\ref{eq:feasCut}) is straightforward. Regarding (\ref{eq:optCut}), a partial solution that replicates the input to the subproblem from iteration $n$ (i.e., maintaining the same paths for all requests) results in an objective value $z$ refined to the corresponding upper bound $\bar{\zeta}^{n}$. Any deviation in the paths neutralises the effect of cut (\ref{eq:optCut}), thus ensuring its validity. Similarly, the feasibility cut (\ref{eq:feasCut}) guarantees that at least one edge from any previously retrieved set of paths is eliminated. We add cuts at any node in which an integer solution is found and maintain cuts in all subsequent nodes.
\section{Large Neighbourhood Search} \label{sec:lns}
    
We propose a metaheuristic methodology that builds upon and enhances existing techniques. Our aim is to develop a simple yet effective method, ensuring that solution quality is not heavily dependent on the specific characteristics of each problem instance. Following this aim, we propose a method based on LNS that employs a single combination of destroy and repair operators. Unlike existing ALNS techniques, our approach eliminates the need for the parameter tuning associated with the adaptive selection of destroy and repair mechanisms. Furthermore, we introduce targeted modifications to widely-used destroy and repair operators to enhance the robustness of our method. Our computational experiences appears to support our adaptation (see Section \ref{sec:experiments}).

\subsection{Destroy operator}

\textbf{Degree of destruction. } The degree of destruction applied to the solution is a crucial LNS parameter \cite{Pisinger2019-wh}. If it is too small, the neighborhood defined by the destroy operator becomes narrow, restricting the exploration of the search space. Conversely, if it is too large, the local search process may resemble a re-optimisation procedure rather than an effective exploration strategy. To determine this, we adopt the approach of \cite{Ropke2006-wm}, where the number of removed requests is randomly selected within a predefined range.

\noindent\textbf{Request selection criterion. }For the selection criterion of the requests to be removed, we utilise a modified version of the related removal operator proposed in \cite{Danloup2018-lh}. The underlying principle of this operator is to remove requests that exhibit low dissimilarity, as this increases the likelihood of finding a beneficial reinsertion that improves the overall solution cost \cite{Shaw1998-xk}.
        

Let $q_r$ be the demand associated with request $r$. Let $p_r$ and $d_r$ be the pickup and delivery node associated with $r$, respectively. Let $l_{i}$ be the earliest possible service time for a node $i$. Let $t_{ij}$ be the travel time associated with traveling from $i$ to $j$. The dissimilarity between requests $r$ and $r'$ in the related removal operator is defined as:
        $$ D(r, r') = \theta_1 \cdot |q_r-q_{r'}| + \theta_2 \cdot \left( t_{p_{r}p_{r'}} + t_{d_{r}d_{r'}} \right) +  \theta_3 \cdot \left( |l_{p_r} - l_{p_{r'}}| + |l_{d_r} - l_{d_{r'}}| \right),$$
        where $\theta_1$, $\theta_2$ and $\theta_3$ are parameters to be tuned, primarily used to handle the potentially different scales of demands, travel times and time windows.
        
 Our proposed modification refines the computation of request dissimilarity, eliminating the need for parameter tuning and thereby simplifying the algorithm. Specifically, each request $r$ is represented as a feature vector
        \[vec(r) =
        \begin{bmatrix}
            q_r & \textsc{x}(p_r) & \textsc{y}(p_r) & \textsc{x}(d_r) & \textsc{y}(d_r) & l_{p_r} & l_{d_r} & st(p_r) & st(d_r)
        \end{bmatrix}^T,\]
        where $\textsc{x}(i)$ and $\textsc{y}(i)$ are the coordinates of node $i$ respectively, and $st(i)$ is its service time.
        
Dissimilarity between requests is measured using the Mahalanobis distance between their vector representations. Mahalanobis distance eliminates the influence of scale, while also accounting for correlations between features, preventing the over-weighting of redundant information.

\subsection{Repair operator}

\textbf{Insertion order. } To determine the order in which the non-serviced requests are to be considered, we use a modified version of the `Insertion Ease' metric proposed by \cite{Danloup2018-lh}. The key idea is that requests that are hard with respect to the problem's constraints, should be inserted in the partial solution first. The Insertion Ease of a request $r$ is defined as follows:
    $$    IE(r)=\theta_1 \cdot q_r + \theta_2 \cdot t_{p_{r}d_{r}} - \theta_3 \cdot (u_{p_r} - l_{p_r} + u_{d_r} - l_{d_r}),
       $$
\ where $u_{i}$ is the latest possible time to begin service at node $i$, and $\theta_1$, $\theta_2$ and $\theta_3$ are parameters to be tuned, primarily used to handle the potentially different scales of demands, travel times and time windows.
        
To enhance the robustness of this approach and avoid problem-specific parameter tuning, we propose a modified formula to compute the insertion difficulty of each request. Specifically, we first compute the minimum and maximum of the following quantities among the instance's requests: demand $q_r$, pickup time window width $w(p_r)$, delivery time window width $w(d_r)$, pickup service time $st(p_r)$, delivery service time $st(d_r)$, and travel time between the pickup and delivery node $t_{p_{r}d_{r}}$. We define the Insertion Difficulty of request $r$ as follows:
        \[ID(r) = \overline{q}_r + \overline{t}_{p_{r}d_{r}}
          + \overline{st}(p_r) + \overline{st}(d_r)
          - \overline{w}(p_r) - \overline{w}(d_r),\]
        where $\overline{x}$ denotes that min-max scaling was applied to the quantity $x$. We sort the unserviced requests in a decreasing Insertion Difficulty order.
        
\noindent\textbf{Request insertion. } For the insertion of each request $r$ in the partial solution $s^-$, we use a randomized cheapest-insertion heuristic. That is, every possible insertion of $r$ is evaluated, with or without the use of a transfer point $t$, leading to a complexity of $O(|\mathcal{T}| \cdot |\mathcal{K}|^2 \cdot |\mathcal{R}|^4)$, which is far from negligible. Checking the feasibility of each candidate insertion is done through the method proposed by \cite{Masson2013-xm}, i.e., using constant-time feasibility testing. Now, to introduce randomness to the insertion procedure, which has proven to be advantageous for the optimisation process \cite{Pisinger2019-wh}, we employ the technique proposed by \cite{Christiaens2020-vb}. Specifically, each feasible candidate insertion is omitted with a probability $\beta$. Thus, the insertion heuristic returns the cheapest feasible insertion not-skipped by the randomisation process. 

\subsection{Acceptance criterion}
Allowing the consideration of non-improving solutions has been widely recognized for its advantages and has become an essential component of many local search techniques. 
Although Simulated Annealing \cite{Kirkpatrick1983-jb} and its variants are relevant acceptance mechanisms, they require extensive parameter tuning that might become highly problem-specific and hence significantly impact performance. We therefore adopt an alternative acceptance criterion of \cite{Burke2017-vf} to minimise tuning and, hopefully, remain robust across various instances.
        
Specifically, we accept a new solution if it is superior to the current one or, if worse, provided that it is better than the solution from $L$ iterations prior. By incorporating this delayed acceptance mechanism, we allow temporary degradations in solution quality to escape local optima. The algorithm maintains a historical record of past solutions, iteratively updating it to guide the acceptance of new candidates, thereby improving exploration while requiring minimal parameterisation.

\section{Experiments} \label{sec:experiments}
All experiments run on a server equipped with 32 AMD Ryzen Threadripper PRO 5955WX 16-core processors and 32 GB of RAM, running Ubuntu 22.04.5 LTS. To solve the Branch-and-Check algorithm, we employ the \emph{Gurobi 11.0.3} optimiser via the open-source optimisation library \emph{Pyomo 6.8.0} in \emph{Python 3.10.12}. The Large Neighborhood Search algorithm was implemented in C\# using the .NET 8.0 SDK and runtime.
        
\subsection{Literature instances}

\noindent \textbf{Loose Time Windows. }The first group of instances assumes a universal time window of $[0,999]$ for all requests, where the upper bound imposes the use of multiple vehicles for instances exceeding 20 requests. Instances with fewer than 10 requests are not discussed, as both our LBBD and the MILP \cite{Lyu23} find the optimal solution in approximately one second. The tested instances include $|\mathcal{R}| = 10, 12, 15, 20, 25, 30$ requests, with $|\mathcal{K}| = 2, 3$ vehicles and $|\mathcal{T}| = 1, 2, 3$ transfer points. For cases where $|\mathcal{K}| = 3$, the number of transfer points is fixed at $|\mathcal{T}| = 3$. Per each combination, we report averages over 10 randomly generated datasets.

Table \ref{tab:exactAvg} presents a comparative summary of the average performance of our LBBD method against the MILP \cite{Lyu23} within the time limit of one hour. The `MILP \cite{Lyu23}' set of columns reports the average values of the Lower Bound (LB), Upper Bound (UB), optimality Gap (\%), and computational time (Time), but only for instances where a feasible solution was found; the number of such instances is shown in the `Feasible' column. As LBBD finds solutions in more instances than the MILP, `LBBD - comparison' provides the corresponding results for the LBBD method but only for the instances where the MILP found a feasible solution. 
Last, `LBBD - for all instances' presents the results obtained for LBBD over all instances for which it found a feasible solution.

        \begin{table}[tbh]
            \centering
            \caption{Comparison between the MILP \cite{Lyu23} and the LBBD - loose time windows}
            \label{tab:exactAvg}
            \resizebox{\textwidth}{!}{
            \begin{tabular}{|c|c|c|ccccc|cccc|ccccc|}
            \hline
            \multirow{2}{*}{$|\mathcal{R}|$} &
              \multirow{2}{*}{$|\mathcal{K}|$} &
              \multirow{2}{*}{$|\mathcal{T}|$} &
              \multicolumn{5}{c|}{MILP \cite{Lyu23}} &
              \multicolumn{4}{c|}{LBBD - comparison} &
              \multicolumn{5}{c|}{LBBD - for all instances} \\ \cline{4-17} 
             &
               &
               &
              \multicolumn{1}{c|}{LB} &
              \multicolumn{1}{c|}{UB} &
              \multicolumn{1}{c|}{Gap (\%)} &
              \multicolumn{1}{c|}{Time (s)} &
              Feasible &
              \multicolumn{1}{c|}{LB} &
              \multicolumn{1}{c|}{UB} &
              \multicolumn{1}{c|}{Gap (\%)} &
              Time (s) &
              \multicolumn{1}{c|}{LB} &
              \multicolumn{1}{c|}{UB} &
              \multicolumn{1}{c|}{Gap (\%)} &
              \multicolumn{1}{c|}{Time (s)} &
              Feasible \\ \hline
            10 &
              2 &
              1 &
              \multicolumn{1}{c|}{717.27} &
              \multicolumn{1}{c|}{717.27} &
              \multicolumn{1}{c|}{0.00} &
              \multicolumn{1}{c|}{11} &
              10 &
              \multicolumn{1}{c|}{717.27} &
              \multicolumn{1}{c|}{717.27} &
              \multicolumn{1}{c|}{0.00} &
              3 &
              \multicolumn{1}{c|}{717.27} &
              \multicolumn{1}{c|}{717.27} &
              \multicolumn{1}{c|}{0.00} &
              \multicolumn{1}{c|}{3} &
              10 \\ \hline
            10 &
              2 &
              2 &
              \multicolumn{1}{c|}{703.04} &
              \multicolumn{1}{c|}{703.04} &
              \multicolumn{1}{c|}{0.00} &
              \multicolumn{1}{c|}{21} &
              10 &
              \multicolumn{1}{c|}{703.04} &
              \multicolumn{1}{c|}{703.04} &
              \multicolumn{1}{c|}{0.00} &
              10 &
              \multicolumn{1}{c|}{703.04} &
              \multicolumn{1}{c|}{703.04} &
              \multicolumn{1}{c|}{0.00} &
              \multicolumn{1}{c|}{9} &
              10 \\ \hline
            10 &
              3 &
              3 &
              \multicolumn{1}{c|}{686.75} &
              \multicolumn{1}{c|}{686.75} &
              \multicolumn{1}{c|}{0.00} &
              \multicolumn{1}{c|}{8} &
              10 &
              \multicolumn{1}{c|}{686.75} &
              \multicolumn{1}{c|}{686.75} &
              \multicolumn{1}{c|}{0.00} &
              26 &
              \multicolumn{1}{c|}{686.75} &
              \multicolumn{1}{c|}{686.75} &
              \multicolumn{1}{c|}{0.00} &
              \multicolumn{1}{c|}{25} &
              10 \\ \hline
            12 &
              2 &
              1 &
              \multicolumn{1}{c|}{729.06} &
              \multicolumn{1}{c|}{729.06} &
              \multicolumn{1}{c|}{0.00} &
              \multicolumn{1}{c|}{43} &
              10 &
              \multicolumn{1}{c|}{729.06} &
              \multicolumn{1}{c|}{729.06} &
              \multicolumn{1}{c|}{0.00} &
              3 &
              \multicolumn{1}{c|}{729.06} &
              \multicolumn{1}{c|}{729.06} &
              \multicolumn{1}{c|}{0.00} &
              \multicolumn{1}{c|}{3} &
              10 \\ \hline
            12 &
              2 &
              2 &
              \multicolumn{1}{c|}{751.60} &
              \multicolumn{1}{c|}{751.60} &
              \multicolumn{1}{c|}{0.00} &
              \multicolumn{1}{c|}{72} &
              10 &
              \multicolumn{1}{c|}{751.60} &
              \multicolumn{1}{c|}{751.60} &
              \multicolumn{1}{c|}{0.00} &
              11 &
              \multicolumn{1}{c|}{751.60} &
              \multicolumn{1}{c|}{751.60} &
              \multicolumn{1}{c|}{0.00} &
              \multicolumn{1}{c|}{11} &
              10 \\ \hline
            12 &
              3 &
              3 &
              \multicolumn{1}{c|}{812.36} &
              \multicolumn{1}{c|}{812.36} &
              \multicolumn{1}{c|}{0.00} &
              \multicolumn{1}{c|}{218} &
              10 &
              \multicolumn{1}{c|}{812.36} &
              \multicolumn{1}{c|}{812.36} &
              \multicolumn{1}{c|}{0.00} &
              154 &
              \multicolumn{1}{c|}{812.36} &
              \multicolumn{1}{c|}{812.36} &
              \multicolumn{1}{c|}{0.00} &
              \multicolumn{1}{c|}{154} &
              10 \\ \hline
            15 &
              2 &
              1 &
              \multicolumn{1}{c|}{933.01} &
              \multicolumn{1}{c|}{933.01} &
              \multicolumn{1}{c|}{0.00} &
              \multicolumn{1}{c|}{362} &
              10 &
              \multicolumn{1}{c|}{933.01} &
              \multicolumn{1}{c|}{933.01} &
              \multicolumn{1}{c|}{0.00} &
              13 &
              \multicolumn{1}{c|}{933.01} &
              \multicolumn{1}{c|}{933.01} &
              \multicolumn{1}{c|}{0.00} &
              \multicolumn{1}{c|}{13} &
              10 \\ \hline
            15 &
              2 &
              2 &
              \multicolumn{1}{c|}{923.66} &
              \multicolumn{1}{c|}{923.66} &
              \multicolumn{1}{c|}{0.00} &
              \multicolumn{1}{c|}{108} &
              10 &
              \multicolumn{1}{c|}{923.66} &
              \multicolumn{1}{c|}{923.66} &
              \multicolumn{1}{c|}{0.00} &
              24 &
              \multicolumn{1}{c|}{923.66} &
              \multicolumn{1}{c|}{923.66} &
              \multicolumn{1}{c|}{0.00} &
              \multicolumn{1}{c|}{23} &
              10 \\ \hline
            15 &
              3 &
              3 &
              \multicolumn{1}{c|}{946.38} &
              \multicolumn{1}{c|}{948.08} &
              \multicolumn{1}{c|}{0.19} &
              \multicolumn{1}{c|}{589} &
              10 &
              \multicolumn{1}{c|}{948.08} &
              \multicolumn{1}{c|}{948.08} &
              \multicolumn{1}{c|}{0.00} &
              442 &
              \multicolumn{1}{c|}{948.08} &
              \multicolumn{1}{c|}{948.08} &
              \multicolumn{1}{c|}{0.00} &
              \multicolumn{1}{c|}{442} &
              10 \\ \hline
            20 &
              2 &
              1 &
              \multicolumn{1}{c|}{1135.55} &
              \multicolumn{1}{c|}{1163.69} &
              \multicolumn{1}{c|}{2.81} &
              \multicolumn{1}{c|}{1985} &
              9 &
              \multicolumn{1}{c|}{1154.80} &
              \multicolumn{1}{c|}{1154.80} &
              \multicolumn{1}{c|}{0.00} &
              143 &
              \multicolumn{1}{c|}{1157.03} &
              \multicolumn{1}{c|}{1157.03} &
              \multicolumn{1}{c|}{0.00} &
              \multicolumn{1}{c|}{178} &
              10 \\ \hline
            20 &
              2 &
              2 &
              \multicolumn{1}{c|}{1144.24} &
              \multicolumn{1}{c|}{1155.54} &
              \multicolumn{1}{c|}{0.99} &
              \multicolumn{1}{c|}{2292} &
              6 &
              \multicolumn{1}{c|}{1154.56} &
              \multicolumn{1}{c|}{1154.56} &
              \multicolumn{1}{c|}{0.00} &
              271 &
              \multicolumn{1}{c|}{1116.39} &
              \multicolumn{1}{c|}{1128.70} &
              \multicolumn{1}{c|}{1.10} &
              \multicolumn{1}{c|}{969} &
              10 \\ \hline
            20 &
              3 &
              3 &
              \multicolumn{1}{c|}{1083.58} &
              \multicolumn{1}{c|}{1186.04} &
              \multicolumn{1}{c|}{8.14} &
              \multicolumn{1}{c|}{2616} &
              9 &
              \multicolumn{1}{c|}{1113.37} &
              \multicolumn{1}{c|}{1114.55} &
              \multicolumn{1}{c|}{0.12} &
              938 &
              \multicolumn{1}{c|}{1103.19} &
              \multicolumn{1}{c|}{1105.76} &
              \multicolumn{1}{c|}{0.25} &
              \multicolumn{1}{c|}{1204} &
              10 \\ \hline
            25 &
              2 &
              1 &
              \multicolumn{1}{c|}{1381.81} &
              \multicolumn{1}{c|}{1439.50} &
              \multicolumn{1}{c|}{4.26} &
              \multicolumn{1}{c|}{2393} &
              2 &
              \multicolumn{1}{c|}{1416.85} &
              \multicolumn{1}{c|}{1416.85} &
              \multicolumn{1}{c|}{0.00} &
              40 &
              \multicolumn{1}{c|}{1291.99} &
              \multicolumn{1}{c|}{1329.76} &
              \multicolumn{1}{c|}{3.29} &
              \multicolumn{1}{c|}{1852} &
              10 \\ \hline
            25 &
              2 &
              2 &
              \multicolumn{1}{c|}{1279.85} &
              \multicolumn{1}{c|}{1395.43} &
              \multicolumn{1}{c|}{8.49} &
              \multicolumn{1}{c|}{3422} &
              2 &
              \multicolumn{1}{c|}{1308.23} &
              \multicolumn{1}{c|}{1329.79} &
              \multicolumn{1}{c|}{1.75} &
              1940 &
              \multicolumn{1}{c|}{1289.43} &
              \multicolumn{1}{c|}{1312.19} &
              \multicolumn{1}{c|}{1.72} &
              \multicolumn{1}{c|}{2572} &
              10 \\ \hline
            25 &
              3 &
              3 &
              \multicolumn{1}{c|}{1337.46} &
              \multicolumn{1}{c|}{1506.54} &
              \multicolumn{1}{c|}{11.23} &
              \multicolumn{1}{c|}{3600} &
              3 &
              \multicolumn{1}{c|}{1376.36} &
              \multicolumn{1}{c|}{1376.36} &
              \multicolumn{1}{c|}{0.00} &
              777 &
              \multicolumn{1}{c|}{1311.04} &
              \multicolumn{1}{c|}{1326.76} &
              \multicolumn{1}{c|}{1.26} &
              \multicolumn{1}{c|}{1470} &
              10 \\ \hline
            30 &
              2 &
              1 &
              \multicolumn{1}{c|}{1412.27} &
              \multicolumn{1}{c|}{1574.83} &
              \multicolumn{1}{c|}{10.32} &
              \multicolumn{1}{c|}{3600} &
              1 &
              \multicolumn{1}{c|}{1454.08} &
              \multicolumn{1}{c|}{1463.82} &
              \multicolumn{1}{c|}{0.67} &
              3600 &
              \multicolumn{1}{c|}{1534.16} &
              \multicolumn{1}{c|}{1596.54} &
              \multicolumn{1}{c|}{3.99} &
              \multicolumn{1}{c|}{2985} &
              10 \\ \hline
            30 &
              2 &
              2 &
              \multicolumn{1}{c|}{-} &
              \multicolumn{1}{c|}{-} &
              \multicolumn{1}{c|}{-} &
              \multicolumn{1}{c|}{-} &
              0 &
              \multicolumn{1}{c|}{-} &
              \multicolumn{1}{c|}{-} &
              \multicolumn{1}{c|}{-} &
              - &
              \multicolumn{1}{c|}{1425.03} &
              \multicolumn{1}{c|}{1508.49} &
              \multicolumn{1}{c|}{5.71} &
              \multicolumn{1}{c|}{3376} &
              7 \\ \hline
            30 &
              3 &
              3 &
              \multicolumn{1}{c|}{-} &
              \multicolumn{1}{c|}{-} &
              \multicolumn{1}{c|}{-} &
              \multicolumn{1}{c|}{-} &
              0 &
              \multicolumn{1}{c|}{-} &
              \multicolumn{1}{c|}{-} &
              \multicolumn{1}{c|}{-} &
              - &
              \multicolumn{1}{c|}{1479.26} &
              \multicolumn{1}{c|}{1560.68} &
              \multicolumn{1}{c|}{5.53} &
              \multicolumn{1}{c|}{2943} &
              9 \\ \hline
            \end{tabular}
            }
        \end{table}

For instances with 10 and 12 requests, both methods solve all datasets to optimality within seconds, making the comparison inconclusive. Similarly, for 29 out of 30 instances with 15 requests, both approaches yield optimal solutions. The only exception is a dataset with 3 vehicles and 3 transfer points, which the MILP fails to solve optimally - hence, the average optimality gap for the MILP is slightly above zero. 

Nevertheless, for larger instances, the LBBD consistently outperforms the MILP in both lower and upper bounds. As a result, the LBBD maintains a notably smaller average optimality gap, while the MILP’s performance deteriorates. Additionally, the LBBD produces feasible solutions in 16 out of 18 instance groups at an average gap of at most 4\%. In the remaining two groups, a few instances remain unsolved for LBBD within one hour, whereas MILP solves none. Detailed results are displayed in Table \ref{tab:exactAll} of the Appendix.

Next, for the larger datasets of \cite{Lyu23} that include 20-30 requests, the LNS algorithm is applied first to generate a feasible solution. This solution is then used as a warm-start for the master problem in the Branch-and-Check, allowing LBBD to explore a branching tree that originates from a feasible node. Table \ref{tab:lbbdWS} presents relevant results for three approaches: the Branch-and-Check algorithm (`LBBD'), the LNS algorithm (`LNS'), and the Branch-and-Check algorithm initialised with a warm-start solution from LNS (`LBBD$^+$'). Each row reports average values across 10 instances of the same number of requests, vehicles, and transfer points. It is worth noting that for four specific instances - three with 30 requests, 2 vehicles, and 2 transfer points, and one with 30 requests, 3 vehicles, and 3 transfer points - the standalone LBBD algorithm failed to find a feasible solution within the time limit. Consequently, these instances are excluded from the average values in the `LNS' and `LBBD$^+$' columns to enable comparison, despite both methods successfully identifying feasible solutions. Table \ref{tab:detLBBDws} of the Appendix presents the results in detail.

For the LNS results, the column `Err' refers to the average deviation between the upper bounds of LNS and LBBD, computed as $\text{Err} = \frac{\text{UB(LNS) - UB(LBBD)}}{\text{UB(LNS)}}\%$. Column `Gap' reports the optimality gap, defined as the relative difference between the upper bound of LNS and the best lower bound, as provided consistently by `LBBD$^+$'. These metrics demonstrate that LNS produces near-optimal solutions in negligible time, thus being an effective warm-start approach. The benefit of warm-starting is further highlighted in columns `Avg Imp' and `Avg Acc'. `Avg Imp' reflects the average improvement in optimality gap for instances that could not be solved to optimality by `LBBD': $\text{Avg Imp = Gap(LBBD) - Gap(LBBD$^+$)}\%$. Meanwhile, `Avg Acc' captures the average time savings for instances successfully solved to optimality by both methods: $\text{Avg Acc = Time(LBBD)}$ $\text{ - Time(LBBD$^+$)}$. Note that `Avg Imp' is not reported for instance groups where all 10 cases are solved to optimality (i.e., for 20 requests, 2 vehicles and 1 transfer point). Overall, the results underscore the value of integrating LNS as a warm-start mechanism within the LBBD framework. 
\begin{table}[tbh]
\centering
\caption{Average results of LNS and LBBD with or without a warm-start solution}
\label{tab:lbbdWS}
\resizebox{\textwidth}{!}{
\begin{tabular}{|c|c|c|cccc|cccc|cccc|c|c|}
\hline
\multirow{2}{*}{$|\mathcal{R}|$} &
  \multirow{2}{*}{$|\mathcal{K}|$} &
  \multirow{2}{*}{$|\mathcal{T}|$} &
  \multicolumn{4}{c|}{LBBD} &
  \multicolumn{4}{c|}{LNS} &
  \multicolumn{4}{c|}{LBBD$^+$} &
  \multirow{2}{*}{Avg Imp} &
  \multirow{2}{*}{Avg Acc} \\ \cline{4-15}
 &
   &
   &
  \multicolumn{1}{c|}{LB} &
  \multicolumn{1}{c|}{UB} &
  \multicolumn{1}{c|}{Gap} &
  Time &
  \multicolumn{1}{c|}{UB} &
  \multicolumn{1}{c|}{Err} &
  \multicolumn{1}{c|}{Gap} &
  Time &
  \multicolumn{1}{c|}{LB} &
  \multicolumn{1}{c|}{UB} &
  \multicolumn{1}{c|}{Gap} &
  Time &
   &
   \\ \hline
20 &
  2 &
  1 &
  \multicolumn{1}{c|}{1157.03} &
  \multicolumn{1}{c|}{1157.03} &
  \multicolumn{1}{c|}{0.00} &
  179 &
  \multicolumn{1}{c|}{1199.12} &
  \multicolumn{1}{c|}{3.47} &
  \multicolumn{1}{c|}{3.47} &
  1 &
  \multicolumn{1}{c|}{1157.03} &
  \multicolumn{1}{c|}{1157.03} &
  \multicolumn{1}{c|}{0.00} &
  129 &
  - &
  50 \\ \hline
20 &
  2 &
  2 &
  \multicolumn{1}{c|}{1116.39} &
  \multicolumn{1}{c|}{1128.70} &
  \multicolumn{1}{c|}{1.10} &
  969 &
  \multicolumn{1}{c|}{1173.72} &
  \multicolumn{1}{c|}{3.81} &
  \multicolumn{1}{c|}{4.83} &
  1 &
  \multicolumn{1}{c|}{1117.05} &
  \multicolumn{1}{c|}{1128.26} &
  \multicolumn{1}{c|}{1.00} &
  907 &
  0.49 &
  77 \\ \hline
20 &
  3 &
  3 &
  \multicolumn{1}{c|}{1103.19} &
  \multicolumn{1}{c|}{1105.76} &
  \multicolumn{1}{c|}{0.25} &
  1204 &
  \multicolumn{1}{c|}{1156.84} &
  \multicolumn{1}{c|}{4.40} &
  \multicolumn{1}{c|}{4.63} &
  1 &
  \multicolumn{1}{c|}{1103.37} &
  \multicolumn{1}{c|}{1105.76} &
  \multicolumn{1}{c|}{0.24} &
  1144 &
  0.09 &
  76 \\ \hline
25 &
  2 &
  1 &
  \multicolumn{1}{c|}{1291.99} &
  \multicolumn{1}{c|}{1329.76} &
  \multicolumn{1}{c|}{3.29} &
  1852 &
  \multicolumn{1}{c|}{1387.14} &
  \multicolumn{1}{c|}{3.92} &
  \multicolumn{1}{c|}{6.61} &
  1 &
  \multicolumn{1}{c|}{1297.97} &
  \multicolumn{1}{c|}{1323.58} &
  \multicolumn{1}{c|}{2.24} &
  1838 &
  2.09 &
  28 \\ \hline
25 &
  2 &
  2 &
  \multicolumn{1}{c|}{1289.43} &
  \multicolumn{1}{c|}{1312.19} &
  \multicolumn{1}{c|}{1.72} &
  2573 &
  \multicolumn{1}{c|}{1367.67} &
  \multicolumn{1}{c|}{4.02} &
  \multicolumn{1}{c|}{5.39} &
  1 &
  \multicolumn{1}{c|}{1293.43} &
  \multicolumn{1}{c|}{1312.19} &
  \multicolumn{1}{c|}{1.42} &
  2429 &
  0.51 &
  361 \\ \hline
25 &
  3 &
  3 &
  \multicolumn{1}{c|}{1311.04} &
  \multicolumn{1}{c|}{1326.76} &
  \multicolumn{1}{c|}{1.26} &
  1470 &
  \multicolumn{1}{c|}{1391.22} &
  \multicolumn{1}{c|}{4.65} &
  \multicolumn{1}{c|}{5.73} &
  1 &
  \multicolumn{1}{c|}{1312.66} &
  \multicolumn{1}{c|}{1326.05} &
  \multicolumn{1}{c|}{1.08} &
  1401 &
  0.60 &
  98 \\ \hline
30 &
  2 &
  1 &
  \multicolumn{1}{c|}{1534.16} &
  \multicolumn{1}{c|}{1596.54} &
  \multicolumn{1}{c|}{3.99} &
  2985 &
  \multicolumn{1}{c|}{1692.71} &
  \multicolumn{1}{c|}{5.65} &
  \multicolumn{1}{c|}{9.15} &
  1 &
  \multicolumn{1}{c|}{1539.56} &
  \multicolumn{1}{c|}{1592.38} &
  \multicolumn{1}{c|}{3.40} &
  2888 &
  0.85 &
  323 \\ \hline
30 &
  2 &
  2 &
  \multicolumn{1}{c|}{1425.03} &
  \multicolumn{1}{c|}{1508.49} &
  \multicolumn{1}{c|}{5.71} &
  3376 &
  \multicolumn{1}{c|}{1602.89} &
  \multicolumn{1}{c|}{5.96} &
  \multicolumn{1}{c|}{10.89} &
  1 &
  \multicolumn{1}{c|}{1431.97} &
  \multicolumn{1}{c|}{1493.79} &
  \multicolumn{1}{c|}{4.32} &
  3316 &
  1.63 &
  422 \\ \hline
30 &
  3 &
  3 &
  \multicolumn{1}{c|}{1479.26} &
  \multicolumn{1}{c|}{1560.68} &
  \multicolumn{1}{c|}{5.53} &
  2944 &
  \multicolumn{1}{c|}{1657.13} &
  \multicolumn{1}{c|}{5.55} &
  \multicolumn{1}{c|}{10.62} &
  1 &
  \multicolumn{1}{c|}{1484.55} &
  \multicolumn{1}{c|}{1551.97} &
  \multicolumn{1}{c|}{4.64} &
  2925 &
  1.15 &
  197 \\ \hline
\end{tabular}
}
\end{table}

\noindent \textbf{Strict Time Windows. } We consider three different time window spans of 180, 240, and 300 minutes. Within each span, three request generation scenarios are examined: `L' (long distances between pickup and delivery points), `S' (short distances), and `M' (a mix of both). Instance generation follows the methodology of \cite{Sam20}, with a small number of requests (3, 4, and 5 requests) and a fleet of 4 vehicles. Additionally, all datasets include 4 transfer points, while instances with 3 requests are also tested with 5 transfer points.

\begin{table}[tbh]
            \centering
            \caption{Comparison between the MILP \cite{Lyu23} and the LBBD - tight time windows}
            \label{tab:twAvg}
            \resizebox{\textwidth}{!}{
            \begin{tabular}{|c|c|c|ccccc|ccccc|}
            \hline
            \multirow{2}{*}{$|\mathcal{R}|$} &
              \multirow{2}{*}{TW} &
              \multirow{2}{*}{$|\mathcal{T}|$} &
              \multicolumn{5}{c|}{MILP of \cite{Lyu23}} &
              \multicolumn{5}{c|}{LBBD} \\ \cline{4-13} 
             &
               &
               &
              \multicolumn{1}{c|}{LB} &
              \multicolumn{1}{c|}{UB} &
              \multicolumn{1}{c|}{Gap (\%)} &
              \multicolumn{1}{c|}{Time (s)} &
              Feasible &
              \multicolumn{1}{c|}{LB} &
              \multicolumn{1}{c|}{UB} &
              \multicolumn{1}{c|}{Gap (\%)} &
              \multicolumn{1}{c|}{Time (s)} &
              Feasible \\ \hline
            3 &
              180L &
              4 &
              \multicolumn{1}{c|}{411.75} &
              \multicolumn{1}{c|}{411.75} &
              \multicolumn{1}{c|}{0.00} &
              \multicolumn{1}{c|}{7} &
              10 &
              \multicolumn{1}{c|}{411.75} &
              \multicolumn{1}{c|}{411.75} &
              \multicolumn{1}{c|}{0.00} &
              \multicolumn{1}{c|}{23} &
              10 \\ \hline
            3 &
              180M &
              4 &
              \multicolumn{1}{c|}{330.92} &
              \multicolumn{1}{c|}{330.92} &
              \multicolumn{1}{c|}{0.00} &
              \multicolumn{1}{c|}{8} &
              10 &
              \multicolumn{1}{c|}{330.92} &
              \multicolumn{1}{c|}{330.92} &
              \multicolumn{1}{c|}{0.00} &
              \multicolumn{1}{c|}{18} &
              10 \\ \hline
            3 &
              180S &
              4 &
              \multicolumn{1}{c|}{261.33} &
              \multicolumn{1}{c|}{261.33} &
              \multicolumn{1}{c|}{0.00} &
              \multicolumn{1}{c|}{3} &
              10 &
              \multicolumn{1}{c|}{261.33} &
              \multicolumn{1}{c|}{261.33} &
              \multicolumn{1}{c|}{0.00} &
              \multicolumn{1}{c|}{7} &
              10 \\ \hline
            3 &
              240L &
              4 &
              \multicolumn{1}{c|}{376.67} &
              \multicolumn{1}{c|}{376.67} &
              \multicolumn{1}{c|}{0.00} &
              \multicolumn{1}{c|}{4} &
              10 &
              \multicolumn{1}{c|}{376.67} &
              \multicolumn{1}{c|}{376.67} &
              \multicolumn{1}{c|}{0.00} &
              \multicolumn{1}{c|}{10} &
              10 \\ \hline
            3 &
              240M &
              4 &
              \multicolumn{1}{c|}{299.72} &
              \multicolumn{1}{c|}{299.72} &
              \multicolumn{1}{c|}{0.00} &
              \multicolumn{1}{c|}{4} &
              10 &
              \multicolumn{1}{c|}{299.72} &
              \multicolumn{1}{c|}{299.72} &
              \multicolumn{1}{c|}{0.00} &
              \multicolumn{1}{c|}{7} &
              10 \\ \hline
            3 &
              240S &
              4 &
              \multicolumn{1}{c|}{251.64} &
              \multicolumn{1}{c|}{251.64} &
              \multicolumn{1}{c|}{0.00} &
              \multicolumn{1}{c|}{2} &
              10 &
              \multicolumn{1}{c|}{251.64} &
              \multicolumn{1}{c|}{251.64} &
              \multicolumn{1}{c|}{0.00} &
              \multicolumn{1}{c|}{2} &
              10 \\ \hline
            3 &
              300L &
              4 &
              \multicolumn{1}{c|}{376.45} &
              \multicolumn{1}{c|}{376.45} &
              \multicolumn{1}{c|}{0.00} &
              \multicolumn{1}{c|}{6} &
              10 &
              \multicolumn{1}{c|}{376.45} &
              \multicolumn{1}{c|}{376.45} &
              \multicolumn{1}{c|}{0.00} &
              \multicolumn{1}{c|}{9} &
              10 \\ \hline
            3 &
              300M &
              4 &
              \multicolumn{1}{c|}{299.72} &
              \multicolumn{1}{c|}{299.72} &
              \multicolumn{1}{c|}{0.00} &
              \multicolumn{1}{c|}{4} &
              10 &
              \multicolumn{1}{c|}{299.72} &
              \multicolumn{1}{c|}{299.72} &
              \multicolumn{1}{c|}{0.00} &
              \multicolumn{1}{c|}{7} &
              10 \\ \hline
            3 &
              300S &
              4 &
              \multicolumn{1}{c|}{251.64} &
              \multicolumn{1}{c|}{251.64} &
              \multicolumn{1}{c|}{0.00} &
              \multicolumn{1}{c|}{2} &
              10 &
              \multicolumn{1}{c|}{251.64} &
              \multicolumn{1}{c|}{251.64} &
              \multicolumn{1}{c|}{0.00} &
              \multicolumn{1}{c|}{2} &
              10 \\ \hline
            3 &
              180L &
              5 &
              \multicolumn{1}{c|}{400.83} &
              \multicolumn{1}{c|}{400.83} &
              \multicolumn{1}{c|}{0.00} &
              \multicolumn{1}{c|}{33} &
              10 &
              \multicolumn{1}{c|}{400.83} &
              \multicolumn{1}{c|}{400.83} &
              \multicolumn{1}{c|}{0.00} &
              \multicolumn{1}{c|}{38} &
              10 \\ \hline
            3 &
              180M &
              5 &
              \multicolumn{1}{c|}{328.53} &
              \multicolumn{1}{c|}{328.53} &
              \multicolumn{1}{c|}{0.00} &
              \multicolumn{1}{c|}{26} &
              10 &
              \multicolumn{1}{c|}{328.53} &
              \multicolumn{1}{c|}{328.53} &
              \multicolumn{1}{c|}{0.00} &
              \multicolumn{1}{c|}{37} &
              10 \\ \hline
            3 &
              180S &
              5 &
              \multicolumn{1}{c|}{261.00} &
              \multicolumn{1}{c|}{261.00} &
              \multicolumn{1}{c|}{0.00} &
              \multicolumn{1}{c|}{6} &
              10 &
              \multicolumn{1}{c|}{261.00} &
              \multicolumn{1}{c|}{261.00} &
              \multicolumn{1}{c|}{0.00} &
              \multicolumn{1}{c|}{13} &
              10 \\ \hline
            3 &
              240L &
              5 &
              \multicolumn{1}{c|}{372.36} &
              \multicolumn{1}{c|}{372.36} &
              \multicolumn{1}{c|}{0.00} &
              \multicolumn{1}{c|}{14} &
              10 &
              \multicolumn{1}{c|}{372.36} &
              \multicolumn{1}{c|}{372.36} &
              \multicolumn{1}{c|}{0.00} &
              \multicolumn{1}{c|}{13} &
              10 \\ \hline
            3 &
              240M &
              5 &
              \multicolumn{1}{c|}{298.66} &
              \multicolumn{1}{c|}{298.66} &
              \multicolumn{1}{c|}{0.00} &
              \multicolumn{1}{c|}{9} &
              10 &
              \multicolumn{1}{c|}{298.66} &
              \multicolumn{1}{c|}{298.66} &
              \multicolumn{1}{c|}{0.00} &
              \multicolumn{1}{c|}{13} &
              10 \\ \hline
            3 &
              240S &
              5 &
              \multicolumn{1}{c|}{251.64} &
              \multicolumn{1}{c|}{251.64} &
              \multicolumn{1}{c|}{0.00} &
              \multicolumn{1}{c|}{4} &
              10 &
              \multicolumn{1}{c|}{251.64} &
              \multicolumn{1}{c|}{251.64} &
              \multicolumn{1}{c|}{0.00} &
              \multicolumn{1}{c|}{5} &
              10 \\ \hline
            3 &
              300L &
              5 &
              \multicolumn{1}{c|}{372.20} &
              \multicolumn{1}{c|}{372.20} &
              \multicolumn{1}{c|}{0.00} &
              \multicolumn{1}{c|}{13} &
              10 &
              \multicolumn{1}{c|}{372.20} &
              \multicolumn{1}{c|}{372.20} &
              \multicolumn{1}{c|}{0.00} &
              \multicolumn{1}{c|}{12} &
              10 \\ \hline
            3 &
              300M &
              5 &
              \multicolumn{1}{c|}{298.66} &
              \multicolumn{1}{c|}{298.66} &
              \multicolumn{1}{c|}{0.00} &
              \multicolumn{1}{c|}{9} &
              10 &
              \multicolumn{1}{c|}{298.66} &
              \multicolumn{1}{c|}{298.66} &
              \multicolumn{1}{c|}{0.00} &
              \multicolumn{1}{c|}{12} &
              10 \\ \hline
            3 &
              300S &
              5 &
              \multicolumn{1}{c|}{251.64} &
              \multicolumn{1}{c|}{251.64} &
              \multicolumn{1}{c|}{0.00} &
              \multicolumn{1}{c|}{4} &
              10 &
              \multicolumn{1}{c|}{251.64} &
              \multicolumn{1}{c|}{251.64} &
              \multicolumn{1}{c|}{0.00} &
              \multicolumn{1}{c|}{5} &
              10 \\ \hline
            4 &
              180L &
              4 &
              \multicolumn{1}{c|}{491.32} &
              \multicolumn{1}{c|}{491.32} &
              \multicolumn{1}{c|}{0.00} &
              \multicolumn{1}{c|}{167} &
              10 &
              \multicolumn{1}{c|}{491.32} &
              \multicolumn{1}{c|}{491.32} &
              \multicolumn{1}{c|}{0.00} &
              \multicolumn{1}{c|}{313} &
              10 \\ \hline
            4 &
              180M &
              4 &
              \multicolumn{1}{c|}{424.71} &
              \multicolumn{1}{c|}{427.00} &
              \multicolumn{1}{c|}{0.44} &
              \multicolumn{1}{c|}{683} &
              10 &
              \multicolumn{1}{c|}{427.00} &
              \multicolumn{1}{c|}{427.00} &
              \multicolumn{1}{c|}{0.00} &
              \multicolumn{1}{c|}{519} &
              10 \\ \hline
            4 &
              180S &
              4 &
              \multicolumn{1}{c|}{334.53} &
              \multicolumn{1}{c|}{334.53} &
              \multicolumn{1}{c|}{0.00} &
              \multicolumn{1}{c|}{118} &
              10 &
              \multicolumn{1}{c|}{334.53} &
              \multicolumn{1}{c|}{334.53} &
              \multicolumn{1}{c|}{0.00} &
              \multicolumn{1}{c|}{50} &
              10 \\ \hline
            4 &
              240L &
              4 &
              \multicolumn{1}{c|}{445.79} &
              \multicolumn{1}{c|}{445.79} &
              \multicolumn{1}{c|}{0.00} &
              \multicolumn{1}{c|}{85} &
              10 &
              \multicolumn{1}{c|}{445.79} &
              \multicolumn{1}{c|}{445.79} &
              \multicolumn{1}{c|}{0.00} &
              \multicolumn{1}{c|}{68} &
              10 \\ \hline
            4 &
              240M &
              4 &
              \multicolumn{1}{c|}{374.48} &
              \multicolumn{1}{c|}{374.48} &
              \multicolumn{1}{c|}{0.00} &
              \multicolumn{1}{c|}{63} &
              10 &
              \multicolumn{1}{c|}{374.48} &
              \multicolumn{1}{c|}{374.48} &
              \multicolumn{1}{c|}{0.00} &
              \multicolumn{1}{c|}{50} &
              10 \\ \hline
            4 &
              240S &
              4 &
              \multicolumn{1}{c|}{322.41} &
              \multicolumn{1}{c|}{322.41} &
              \multicolumn{1}{c|}{0.00} &
              \multicolumn{1}{c|}{51} &
              10 &
              \multicolumn{1}{c|}{322.41} &
              \multicolumn{1}{c|}{322.41} &
              \multicolumn{1}{c|}{0.00} &
              \multicolumn{1}{c|}{14} &
              10 \\ \hline
            4 &
              300L &
              4 &
              \multicolumn{1}{c|}{444.67} &
              \multicolumn{1}{c|}{444.67} &
              \multicolumn{1}{c|}{0.00} &
              \multicolumn{1}{c|}{88} &
              10 &
              \multicolumn{1}{c|}{444.67} &
              \multicolumn{1}{c|}{444.67} &
              \multicolumn{1}{c|}{0.00} &
              \multicolumn{1}{c|}{65} &
              10 \\ \hline
            4 &
              300M &
              4 &
              \multicolumn{1}{c|}{374.48} &
              \multicolumn{1}{c|}{374.48} &
              \multicolumn{1}{c|}{0.00} &
              \multicolumn{1}{c|}{67} &
              10 &
              \multicolumn{1}{c|}{374.48} &
              \multicolumn{1}{c|}{374.48} &
              \multicolumn{1}{c|}{0.00} &
              \multicolumn{1}{c|}{49} &
              10 \\ \hline
            4 &
              300S &
              4 &
              \multicolumn{1}{c|}{322.41} &
              \multicolumn{1}{c|}{322.41} &
              \multicolumn{1}{c|}{0.00} &
              \multicolumn{1}{c|}{44} &
              10 &
              \multicolumn{1}{c|}{322.41} &
              \multicolumn{1}{c|}{322.41} &
              \multicolumn{1}{c|}{0.00} &
              \multicolumn{1}{c|}{13} &
              10 \\ \hline
            5 &
              180L &
              4 &
              \multicolumn{1}{c|}{537.39} &
              \multicolumn{1}{c|}{564.10} &
              \multicolumn{1}{c|}{5.13} &
              \multicolumn{1}{c|}{1996} &
              5 &
              \multicolumn{1}{c|}{564.10} &
              \multicolumn{1}{c|}{564.10} &
              \multicolumn{1}{c|}{0.00} &
              \multicolumn{1}{c|}{878} &
              7 \\ \hline
            5 &
              180M &
              4 &
              \multicolumn{1}{c|}{458.14} &
              \multicolumn{1}{c|}{467.39} &
              \multicolumn{1}{c|}{1.85} &
              \multicolumn{1}{c|}{1849} &
              6 &
              \multicolumn{1}{c|}{457.42} &
              \multicolumn{1}{c|}{467.39} &
              \multicolumn{1}{c|}{2.13} &
              \multicolumn{1}{c|}{976} &
              7 \\ \hline
            5 &
              180S &
              4 &
              \multicolumn{1}{c|}{404.85} &
              \multicolumn{1}{c|}{404.85} &
              \multicolumn{1}{c|}{0.00} &
              \multicolumn{1}{c|}{676} &
              10 &
              \multicolumn{1}{c|}{404.85} &
              \multicolumn{1}{c|}{404.85} &
              \multicolumn{1}{c|}{0.00} &
              \multicolumn{1}{c|}{123} &
              10 \\ \hline
            5 &
              240L &
              4 &
              \multicolumn{1}{c|}{521.78} &
              \multicolumn{1}{c|}{528.26} &
              \multicolumn{1}{c|}{1.30} &
              \multicolumn{1}{c|}{1265} &
              10 &
              \multicolumn{1}{c|}{528.26} &
              \multicolumn{1}{c|}{528.26} &
              \multicolumn{1}{c|}{0.00} &
              \multicolumn{1}{c|}{213} &
              10 \\ \hline
            5 &
              240M &
              4 &
              \multicolumn{1}{c|}{443.19} &
              \multicolumn{1}{c|}{452.04} &
              \multicolumn{1}{c|}{1.71} &
              \multicolumn{1}{c|}{1341} &
              10 &
              \multicolumn{1}{c|}{452.04} &
              \multicolumn{1}{c|}{452.04} &
              \multicolumn{1}{c|}{0.00} &
              \multicolumn{1}{c|}{417} &
              10 \\ \hline
            5 &
              240S &
              4 &
              \multicolumn{1}{c|}{389.88} &
              \multicolumn{1}{c|}{389.88} &
              \multicolumn{1}{c|}{0.00} &
              \multicolumn{1}{c|}{283} &
              10 &
              \multicolumn{1}{c|}{389.88} &
              \multicolumn{1}{c|}{389.88} &
              \multicolumn{1}{c|}{0.00} &
              \multicolumn{1}{c|}{36} &
              10 \\ \hline
            5 &
              300L &
              4 &
              \multicolumn{1}{c|}{519.60} &
              \multicolumn{1}{c|}{527.92} &
              \multicolumn{1}{c|}{1.62} &
              \multicolumn{1}{c|}{1353} &
              10 &
              \multicolumn{1}{c|}{526.71} &
              \multicolumn{1}{c|}{526.71} &
              \multicolumn{1}{c|}{0.00} &
              \multicolumn{1}{c|}{173} &
              10 \\ \hline
            5 &
              300M &
              4 &
              \multicolumn{1}{c|}{444.18} &
              \multicolumn{1}{c|}{449.88} &
              \multicolumn{1}{c|}{1.02} &
              \multicolumn{1}{c|}{1252} &
              10 &
              \multicolumn{1}{c|}{449.88} &
              \multicolumn{1}{c|}{449.88} &
              \multicolumn{1}{c|}{0.00} &
              \multicolumn{1}{c|}{315} &
              10 \\ \hline
            5 &
              300S &
              4 &
              \multicolumn{1}{c|}{389.65} &
              \multicolumn{1}{c|}{389.65} &
              \multicolumn{1}{c|}{0.00} &
              \multicolumn{1}{c|}{322} &
              10 &
              \multicolumn{1}{c|}{389.65} &
              \multicolumn{1}{c|}{389.65} &
              \multicolumn{1}{c|}{0.00} &
              \multicolumn{1}{c|}{32} &
              10 \\ \hline
            \end{tabular}
            }
            \end{table}

The results appear in Table \ref{tab:twAvg}. Column $|\mathcal{R}|$ denotes the number of requests, while column $|\mathcal{T}|$ indicates the number of transfer points. The time window span (180, 240, or 300 minutes) and distance type (L, M, or S) are indicated under `TW'. The number of vehicles is omitted, as all datasets assume a fleet of 4 vehicles. Each row reports the average results for all datasets matching the specified attributes where both methods found feasible solutions.

In 34 out of 36 instance groups, both methods successfully find feasible solutions for all 10 instances. However, for the groups with 5 requests and the 180L and 180M configurations, neither method could solve all instances within one hour. Specifically, the reported averages for these groups include only 5-7 instances, as shown in columns `Feasible' that indicate the number of instances for which a feasible solution is found.

Both methods achieve optimality for the vast majority of instances. For datasets with 3 and 4 requests, there is no striking difference. However, as the number of requests increases to 5, solving times grow and the LBBD method starts to converge to optimality more quickly. Notably, several instances with 5 requests are not solved to optimality by the MILP, whereas the LBBD reaches optimality for nearly all but the instance groups where neither method could find feasible solutions for all instances. These results suggest that the LBBD is more scalable, making it a strong candidate for solving larger instances efficiently. Table \ref{tab:twAll} in the Appendix, presents the detailed results.

\subsection{Large instances}
\textbf{Generator of benchmarks. }To evaluate the performance of our approach, we develop an instance generator for PDPT. Each generated instance is defined by three sets: a set of transportation requests $\mathcal{R}$, a fleet of vehicles $\mathcal{K}$, and a set of transfer points $\mathcal{T}$. A request $r$ is specified by a pickup location $p_{r}$, a delivery location $d_{r}$, and a cargo quantity $q_{r}$ to be transported between the two. Each vehicle $k$ is defined by its origin $o_{k}$ and destination locations $e_{k}$, along with a maximum cargo capacity $Q$. Transfer points serve as intermediate transfer locations and are defined solely by their geographic coordinates.
        
 Each location in the instance (those associated with requests, vehicles, and transfer points) is characterised by a set of geographical coordinates, a service duration, and a time window indicating the allowable time frame for the start of service. Travel times between locations are computed using the haversine distance (in hecatometers) divided by a fixed vehicle speed, yielding the time in minutes. Regarding the geographical coordinates of the locations, we utilise the street map of Athens, Greece, as provided by OpenStreetMap\footnote{\url{https://www.openstreetmap.org/}}. We retain only the nodes located within a predefined radius from the city's center.
        
The generation of requests follows a feasibility-first approach. Each request is sampled by first proposing pickup and delivery locations, along with associated service durations and time windows. A request is accepted only if there exists a feasible direct route for a single vehicle to depart from its origin location, visit the pickup location and perform the required service, proceed to the delivery location and complete the delivery, and terminate at its destination location. This route must satisfy all relevant time window constraints. If the request passes this check, the cargo quantity is sampled from a predefined distribution, and the request is added to the instance. This process is repeated until $|\mathcal{R}|$ valid requests are generated.
        
To determine the number of vehicles available in the instance, we consider the same problem but without transfers (PDP). We perform a binary search over the number of vehicles to find the smallest fleet size for which a feasible solution can be found through a simple cheapest-insertion heuristic. Finally, we generate transfer locations by performing k-Means clustering on the set of all geographic coordinates associated with the generated requests. The centroids of the resulting clusters are used as the coordinates of the transfer locations.
        
\noindent \textbf{Generated instances. }Using the generator described above, we produce a diverse collection of benchmark instances designed to reflect realistic operational settings and varying levels of complexity. The load $q_{r}$ associated with request $r$ is u.a.r. [5, 25]. The vehicle capacity $Q$ is set to five times the average cargo quantity, resulting in $Q=75$, to ensure that vehicles can serve multiple requests. All time-related quantities are modelled assuming an 8-hour workday, beginning at time 0 and ending at time 480 (minutes). Vehicle travel speed is fixed at 20 km/h. Time is measured in minutes, and travel times are computed by converting haversine distances (in hecatometers) into durations based on this speed.

To generate realistic time windows, the workday is partitioned into 16 half-hour intervals. For each location, the start of its time window is drawn uniformly from the set {0, 30, 60, ..., 450}. The service duration required at each location is u.a.r. [3, 10] minutes. Three levels of time window width were considered: small `S' (60 or 90 minutes), medium `M' (90 or 120 minutes) and large `L' (120 or 150 minutes). Each generated instance uses a fixed time window width selected from one of the above ranges.
        
We generate instances at four different scales, defined by the number of requests and transfer points: 25 requests and 3 transfer points, 50 requests with 4 transfer points, 75 requests with 5 transfer points and 100 requests with 6 transfer points. For each scale and time-window width (small, medium, large), we produce five instances, resulting in a total of 60 instances. This systematic variation enables a thorough analysis of algorithmic performance across diverse problem sizes and time constraints. A summary of the generator's parameters is displayed in Table \ref{tab:generator}. The code of the generator and the instances are publicly available.\footnote{\url{https://github.com/tsompanidisn/PDPTW_T_BENCH}}

        \begin{table}[H]
        \scriptsize
            
            \centering
            \caption{Summary of generator parameters}
            \label{tab:generator}
            \resizebox{\textwidth}{!}{
            \begin{tabular}{rl}
            \hline
            \multicolumn{1}{r}{\textbf{Parameter}} & \multicolumn{1}{l}{\textbf{Uniform}}                              \\
            \hline
            $|\mathcal{R}|$                         & 25, 50, 75, 100                                                   \\
            $q_{r}$                                 & \texttt{integer}[5, 25]                      \\
            $Q$                                     & 75                                                                \\
            $l_{j}$                                 & 0, 30, 60, …, 450 for demand locations, 0 for depots and transfer points                                                 \\
            $u_{j}$                                 & $l_{j}+$: 60 or 90 for `S', 90 or 120 for `M', 120 or 150 for `L' - set to 480 for depots and transfer points \\
            $st(j)$                                 & \texttt{integer}[3, 10]                      \\
            $|\mathcal{T}|$ & 3 for 25 requests, 4 for 50 requests, 5 for 75 requests, 6 for 100 requests \\ \hline
            \end{tabular}}
        \end{table}

\noindent\textbf{Performance of Branch-and-Check. } LBBD is applied to the generated benchmarks of tractable scale for the master problem, i.e., with 25 and 50 requests. Initially, the algorithm is tested as a standalone method. However, due to its difficulty in retrieving a feasible solution independently, a second set of experiments is conducted, where the solution from the refined LNS algorithm is used as a warm-start for the master problem.
        
Table \ref{tab:exact_alns} presents the results for both experimental setups. The column `LBBD' reports the performance of the exact method without warm-start support, while `LBBD$^+$' shows the results after incorporating the LNS solution into the master problem. Additionally, the `rLNS' column provides the objective value of the solution obtained by our LNS method. No time results are given as LBBD reaches the time limit in all instances (with or without warm-start).

        \begin{table}[tbh]
            \centering
            \caption{Results of the Branch-and-Check on the generated benchmarks}
            \label{tab:exact_alns}
            \resizebox{\textwidth}{!}{
            \begin{tabular}{|c|c|c|c|c|ccc|c|ccc|}
            \hline
            \multirow{2}{*}{Requests} &
              \multirow{2}{*}{Vehicles} &
              \multirow{2}{*}{Transfer points} &
              \multirow{2}{*}{Time windows} &
              \multirow{2}{*}{Variant} &
              \multicolumn{3}{c|}{LBBD} &
              rLNS &
              \multicolumn{3}{c|}{LBBD$^+$} \\ \cline{6-12} 
             &
               &
               &
               &
               &
              \multicolumn{1}{c|}{LB} &
              \multicolumn{1}{c|}{UB} &
              \multicolumn{1}{c|}{Gap} &
              UB &
              \multicolumn{1}{c|}{LB} &
              \multicolumn{1}{c|}{UB} &
              \multicolumn{1}{c|}{Gap} \\ \hline
            \multirow{15}{*}{25} &
              3 &
              3 &
              L &
              1 &
              \multicolumn{1}{c|}{992.36} &
              \multicolumn{1}{c|}{-} &
              \multicolumn{1}{c|}{-} &
              1664 &
              \multicolumn{1}{c|}{993.43} &
              \multicolumn{1}{c|}{1664} &
              \multicolumn{1}{c|}{40.30} \\ 
              \cline{2-12} 
             &
              3 &
              3 &
              L &
              2 &
              \multicolumn{1}{c|}{887.94} &
              \multicolumn{1}{c|}{-} &
              \multicolumn{1}{c|}{-} &
              1590 &
              \multicolumn{1}{c|}{\textbf{889.32}} &
              \multicolumn{1}{c|}{\textbf{1423}} &
              \multicolumn{1}{c|}{\textbf{37.50}} 
              \\ \cline{2-12} 
             &
              3 &
              3 &
              L &
              3 &
              \multicolumn{1}{c|}{926.17} &
              \multicolumn{1}{c|}{-} &
              \multicolumn{1}{c|}{-} &
              1585 &
              \multicolumn{1}{c|}{953.96} &
              \multicolumn{1}{c|}{1585} &
              \multicolumn{1}{c|}{39.81}  \\ \cline{2-12} 
             &
              4 &
              3 &
              L &
              4 &
              \multicolumn{1}{c|}{856.43} &
              \multicolumn{1}{c|}{-} &
              \multicolumn{1}{c|}{-} &
              1556 &
              \multicolumn{1}{c|}{\textbf{860.68}} &
              \multicolumn{1}{c|}{\textbf{1526}} &
              \multicolumn{1}{c|}{\textbf{43.60}} \\ \cline{2-12} 
             &
              3 &
              3 &
              L &
              5 &
              \multicolumn{1}{c|}{902.16} &
              \multicolumn{1}{c|}{1610} &
              \multicolumn{1}{c|}{43.97} &
              1404 &
              \multicolumn{1}{c|}{\textbf{904.86}} &
              \multicolumn{1}{c|}{\textbf{1389}} &
              \multicolumn{1}{c|}{\textbf{34.86}} \\ \cline{2-12} 
             &
              4 &
              3 &
              M &
              1 &
              \multicolumn{1}{c|}{1069.95} &
              \multicolumn{1}{c|}{-} &
              \multicolumn{1}{c|}{-} &
              1534 &
              \multicolumn{1}{c|}{\textbf{1070.77}} &
              \multicolumn{1}{c|}{\textbf{1506}} &
              \multicolumn{1}{c|}{\textbf{28.90}}  \\ \cline{2-12} 
             &
              3 &
              3 &
              M &
              2 &
              \multicolumn{1}{c|}{1034.24} &
              \multicolumn{1}{c|}{-} &
              \multicolumn{1}{c|}{-} &
              1674 &
              \multicolumn{1}{c|}{1044.23} &
              \multicolumn{1}{c|}{1674} &
              \multicolumn{1}{c|}{37.62} \\ \cline{2-12} 
             &
              4 &
              3 &
              M &
              3 &
              \multicolumn{1}{c|}{991.66} &
              \multicolumn{1}{c|}{-} &
              \multicolumn{1}{c|}{-} &
              1584 &
              \multicolumn{1}{c|}{1001.00} &
              \multicolumn{1}{c|}{1584} &
              \multicolumn{1}{c|}{36.81} \\ \cline{2-12} 
             &
              4 &
              3 &
              M &
              4 &
              \multicolumn{1}{c|}{899.36} &
              \multicolumn{1}{c|}{-} &
              \multicolumn{1}{c|}{-} &
              1445 &
              \multicolumn{1}{c|}{\textbf{900.21}} &
              \multicolumn{1}{c|}{\textbf{1441}} &
              \multicolumn{1}{c|}{\textbf{37.53}} 
              \\ \cline{2-12} 
             &
              4 &
              3 &
              M &
              5 &
              \multicolumn{1}{c|}{1013.30} &
              \multicolumn{1}{c|}{1781} &
              \multicolumn{1}{c|}{43.10} &
              1855 &
              \multicolumn{1}{c|}{\textbf{1034.71}} &
              \multicolumn{1}{c|}{\textbf{1717}} &
              \multicolumn{1}{c|}{\textbf{39.74}}  \\ \cline{2-12} 
             &
              4 &
              3 &
              S &
              1 &
              \multicolumn{1}{c|}{1058.82} &
              \multicolumn{1}{c|}{1786} &
              \multicolumn{1}{c|}{40.72} &
              1754 &
              \multicolumn{1}{c|}{1082.88} &
              \multicolumn{1}{c|}{1754} &
              \multicolumn{1}{c|}{38.26}  \\ \cline{2-12} 
             &
              4 &
              3 &
              S &
              2 &
              \multicolumn{1}{c|}{1009.01} &
              \multicolumn{1}{c|}{-} &
              \multicolumn{1}{c|}{-} &
              1646 &
              \multicolumn{1}{c|}{1015.89} &
              \multicolumn{1}{c|}{1646} &
              \multicolumn{1}{c|}{38.28}  \\ \cline{2-12} 
             &
              4 &
              3 &
              S &
              3 &
              \multicolumn{1}{c|}{1082.40} &
              \multicolumn{1}{c|}{-} &
              \multicolumn{1}{c|}{-} &
              1649 &
              \multicolumn{1}{c|}{1090.70} &
              \multicolumn{1}{c|}{1649} &
              \multicolumn{1}{c|}{33.86} \\ \cline{2-12} 
             &
              4 &
              3 &
              S &
              4 &
              \multicolumn{1}{c|}{1177.44} &
              \multicolumn{1}{c|}{-} &
              \multicolumn{1}{c|}{-} &
              1682 &
              \multicolumn{1}{c|}{\textbf{1187.45}} &
              \multicolumn{1}{c|}{\textbf{1647}} &
              \multicolumn{1}{c|}{\textbf{27.90}} \\ \cline{2-12} 
             &
              3 &
              3 &
              S &
              5 &
              \multicolumn{1}{c|}{1123.87} &
              \multicolumn{1}{c|}{-} &
              \multicolumn{1}{c|}{-} &
              1757 &
              \multicolumn{1}{c|}{1149.31} &
              \multicolumn{1}{c|}{1757} &
              \multicolumn{1}{c|}{34.59} \\ \hline
            \multirow{15}{*}{50} &
              6 &
              4 &
              L &
              1 &
              \multicolumn{1}{c|}{1441.58} &
              \multicolumn{1}{c|}{-} &
              \multicolumn{1}{c|}{-} &
              2706 &
              \multicolumn{1}{c|}{1442.01} &
              \multicolumn{1}{c|}{2706} &
              \multicolumn{1}{c|}{46.71}  \\ \cline{2-12} 
             &
              6 &
              4 &
              L &
              2 &
              \multicolumn{1}{c|}{1370.09} &
              \multicolumn{1}{c|}{-} &
              \multicolumn{1}{c|}{-} &
              2642 &
              \multicolumn{1}{c|}{1379.90} &
              \multicolumn{1}{c|}{2642} &
              \multicolumn{1}{c|}{47.77} \\ \cline{2-12} 
             &
              5 &
              4 &
              L &
              3 &
              \multicolumn{1}{c|}{1297.64} &
              \multicolumn{1}{c|}{-} &
              \multicolumn{1}{c|}{-} &
              2733 &
              \multicolumn{1}{c|}{1298.51} &
              \multicolumn{1}{c|}{2733} &
              \multicolumn{1}{c|}{52.49} \\ \cline{2-12} 
             &
              7 &
              4 &
              L &
              4 &
              \multicolumn{1}{c|}{1363.51} &
              \multicolumn{1}{c|}{-} &
              \multicolumn{1}{c|}{-} &
              2632 &
              \multicolumn{1}{c|}{1369.77} &
              \multicolumn{1}{c|}{2632} &
              \multicolumn{1}{c|}{47.96} \\ \cline{2-12} 
             &
              6 &
              4 &
              L &
              5 &
              \multicolumn{1}{c|}{1357.78} &
              \multicolumn{1}{c|}{-} &
              \multicolumn{1}{c|}{-} &
              2604 &
              \multicolumn{1}{c|}{1359.99} &
              \multicolumn{1}{c|}{2604} &
              \multicolumn{1}{c|}{47.77} \\ \cline{2-12} 
             &
              7 &
              4 &
              M &
              1 &
              \multicolumn{1}{c|}{1550.29} &
              \multicolumn{1}{c|}{-} &
              \multicolumn{1}{c|}{-} &
              2614 &
              \multicolumn{1}{c|}{1562.00} &
              \multicolumn{1}{c|}{2614} &
              \multicolumn{1}{c|}{40.24} \\ \cline{2-12} 
             &
              7 &
              4 &
              M &
              2 &
              \multicolumn{1}{c|}{1472.47} &
              \multicolumn{1}{c|}{-} &
              \multicolumn{1}{c|}{-} &
              2942 &
              \multicolumn{1}{c|}{1482.12} &
              \multicolumn{1}{c|}{2942} &
              \multicolumn{1}{c|}{49.62} \\ \cline{2-12} 
             &
              7 &
              4 &
              M &
              3 &
              \multicolumn{1}{c|}{1543.21} &
              \multicolumn{1}{c|}{-} &
              \multicolumn{1}{c|}{-} &
              2812 &
              \multicolumn{1}{c|}{\textbf{1550.96}} &
              \multicolumn{1}{c|}{\textbf{2807}} &
              \multicolumn{1}{c|}{\textbf{44.75}} 
              \\ \cline{2-12} 
             &
              7 &
              4 &
              M &
              4 &
              \multicolumn{1}{c|}{1493.33} &
              \multicolumn{1}{c|}{-} &
              \multicolumn{1}{c|}{-} &
              2933 &
              \multicolumn{1}{c|}{\textbf{1493.58}} &
              \multicolumn{1}{c|}{\textbf{2914}} &
              \multicolumn{1}{c|}{\textbf{48.74}} \\ \cline{2-12} 
             &
              7 &
              4 &
              M &
              5 &
              \multicolumn{1}{c|}{1422.21} &
              \multicolumn{1}{c|}{-} &
              \multicolumn{1}{c|}{-} &
              2821 &
              \multicolumn{1}{c|}{1424.77} &
              \multicolumn{1}{c|}{2821} &
              \multicolumn{1}{c|}{49.49} \\ \cline{2-12} 
             &
              7 &
              4 &
              S &
              1 &
              \multicolumn{1}{c|}{1773.74} &
              \multicolumn{1}{c|}{-} &
              \multicolumn{1}{c|}{-} &
              3461 &
              \multicolumn{1}{c|}{\textbf{1786.71}} &
              \multicolumn{1}{c|}{\textbf{3393}} &
              \multicolumn{1}{c|}{\textbf{47.34}} \\ \cline{2-12} 
             &
              6 &
              4 &
              S &
              2 &
              \multicolumn{1}{c|}{1765.88} &
              \multicolumn{1}{c|}{-} &
              \multicolumn{1}{c|}{-} &
              3199 &
              \multicolumn{1}{c|}{\textbf{1780.09}} &
              \multicolumn{1}{c|}{\textbf{3171}} &
              \multicolumn{1}{c|}{\textbf{43.86}} 
              \\ \cline{2-12} 
             &
              7 &
              4 &
              S &
              3 &
              \multicolumn{1}{c|}{1583.91} &
              \multicolumn{1}{c|}{-} &
              \multicolumn{1}{c|}{-} &
              2852 &
              \multicolumn{1}{c|}{\textbf{1586.11}} &
              \multicolumn{1}{c|}{\textbf{2737}} &
              \multicolumn{1}{c|}{\textbf{42.05}}  \\ \cline{2-12} 
             &
              7 &
              4 &
              S &
              4 &
              \multicolumn{1}{c|}{1685.68} &
              \multicolumn{1}{c|}{-} &
              \multicolumn{1}{c|}{-} &
              3196 &
              \multicolumn{1}{c|}{1707.29} &
              \multicolumn{1}{c|}{3196} &
              \multicolumn{1}{c|}{46.58} \\ \cline{2-12} 
             &
              7 &
              4 &
              S &
              5 &
              \multicolumn{1}{c|}{1625.96} &
              \multicolumn{1}{c|}{-} &
              \multicolumn{1}{c|}{-} &
              2947 &
              \multicolumn{1}{c|}{\textbf{1638.86}} &
              \multicolumn{1}{c|}{\textbf{2905}} &
              \multicolumn{1}{c|}{\textbf{43.58}} \\ \hline
            \end{tabular}
            }
\end{table}

The results indicate that the LBBD algorithm struggles to produce feasible solutions, as an upper bound (`UB') is reported for only 3 out of 30 instances, all of which involve 25 requests. When a warm-start solution is provided, a slight improvement in the lower bound (`LB') is observed across all 30 instances. In almost half of the cases (13 out of 30), the LNS solution is improved, yielding a reduced upper bound. The corresponding instances are highlighted in bold.
       
\noindent\textbf{Performance of rLNS. }For larger datasets, we apply the LNS algorithm described in Section \ref{sec:lns}. To assess the effectiveness of rLNS, we conduct a series of computational experiments involving two alternative approaches: `LS' and `MULTI-OP'. The LS metaheuristic is structurally similar to rLNS, with the primary distinction being the use of the original removal operator of \cite{Danloup2018-lh}, in contrast to the modified version used in rLNS. The second baseline, MULTI-OP, is built on the ALNS approach \cite{Ropke2006-wm}, i.e., it operates by iteratively selecting pairs of destroy and repair operators, guided by an adaptive mechanism that updates selection probabilities based on recent performance. The destroy phase of MULTI-OP employs three removal strategies: random removal and worst removal as defined in \cite{Ropke2006-wm}, as well as the similarity-based removal described in Section \ref{sec:lns}. For the repair phase, two request insertion strategies were examined: one involving a random request ordering, and the other using the Insertion Difficulty criterion, as introduced in Section \ref{sec:lns}.

The three methods have been tested on a diverse set of problem instances tailored for PDPT, generated by the aforementioned rules at two different scales, namely 75 and 100 requests. Table~\ref{tab:hyperparameters} summarises the hyperparameter configurations used for each of the three methods evaluated in this experiment. Shared parameters, such as the LAHC \cite{Burke2017-vf} list size and degree of destruction, are kept consistent to ensure a fair comparison, while method-specific parameters reflect the unique mechanisms of each approach. 

\begin{table}[tbh]
        \scriptsize
            
        \centering
        \caption{Hyperparameter settings for each method}
        \label{tab:hyperparameters}
            \begin{tabular}{|l|l|l|}
            \hline
            \textbf{Parameter} & \textbf{Description} & \textbf{Value(s)} \\
            \hline
            \multicolumn{3}{|c|}{\textbf{rLNS}} \\
            \hline
            $L$ & LAHC \cite{Burke2017-vf} list size & 20 \\
            $DR$ & Degree of destruction (uniform range) & [5, 15] \\
            $\beta$ & Blink rate for repair operator & 0.05 \\
            \hline
            \multicolumn{3}{|c|}{\textbf{LS}} \\
            \hline
            $L$ & LAHC \cite{Burke2017-vf} list size & 20 \\
            $DR$ & Degree of destruction (uniform range) & [5, 15] \\
            $\beta$ & Blink rate for repair operator & 0.05 \\
            $\theta_1$, $\theta_2$, $\theta_3$ & Shaw removal operator coefficients & 0.33, 0.99, 0.66 \\
            \hline
            \multicolumn{3}{|c|}{\textbf{MULTI-OP}} \\
            \hline
            $L$ & LAHC \cite{Burke2017-vf} list size & 20 \\
            $DR$ & Degree of destruction (uniform range) & [5, 15] \\
            $\beta$ & Blink rate for repair operator & 0.05 \\
            $w_s$ & Cost scaling factor for worst removal & 0.25 \\
            $\alpha$ & Learning rate for operator weights & 0.3 \\
            $s_{B}$ & Reward for improving best solution & 3 \\
            $s_{A}$ & Reward for accepted solution & 1 \\
            \hline
            \end{tabular}
        \end{table}
Table \ref{tab:meta_comparison} presents the comparative results for MULTI-OP, LS and rLNS. To mitigate the effect of the algorithms' randomness on the experimental results, we follow a multi-restart approach, i.e. each algorithm runs $K = 10$ times on each problem instance. For each restart, the algorithm is terminated when the local search performs $p = 50$ successive iterations without improving the best solution encountered. The `Best UB' columns report the objective value of the best-found solution across all restarts. Columns `Avg UB' indicate the average objective value of the best-found solution per restart, and columns `Time' indicate the average execution time per restart in seconds.

        \begin{table}[tbh]
        \scriptsize
            \centering
            \caption{Comparison between MULTI-OP, LS and the rLNS}
            \label{tab:meta_comparison}
            \resizebox{\textwidth}{!}{
            \begin{tabular}{|c|c|c|ccc|ccc|ccc|}
            \hline
            \multirow{2}{*}{$|\mathcal{R}|$} &
              \multirow{2}{*}{TW} &
              \multirow{2}{*}{Variant} &
              \multicolumn{3}{c|}{MULTI-OP} &
              \multicolumn{3}{c|}{LS} &
              \multicolumn{3}{c|}{rLNS} \\ \cline{4-12} 
             &
               &
               &
              \multicolumn{1}{c|}{Best UB} &
              \multicolumn{1}{c|}{Avg UB} &
              Time &
              \multicolumn{1}{c|}{Best UB} &
              \multicolumn{1}{c|}{Avg UB} &
              Time &
              \multicolumn{1}{c|}{Best UB} &
              \multicolumn{1}{c|}{Avg UB} &
              Time \\ \hline
            \multirow{15}{*}{75} &
              L &
              1 &
              \multicolumn{1}{c|}{3689} &
              \multicolumn{1}{c|}{4038} &
              36 &
              \multicolumn{1}{c|}{3767} &
              \multicolumn{1}{c|}{4051} &
              33 &
              \multicolumn{1}{c|}{\textbf{3644}} &
              \multicolumn{1}{c|}{\textbf{3866}} &
              41 \\ \cline{2-12} 
             &
              L &
              2 &
              \multicolumn{1}{c|}{3429} &
              \multicolumn{1}{c|}{3711} &
              46 &
              \multicolumn{1}{c|}{3590} &
              \multicolumn{1}{c|}{3918} &
              29 &
              \multicolumn{1}{c|}{\textbf{3410}} &
              \multicolumn{1}{c|}{3881} &
              32 \\ \cline{2-12} 
             &
              L &
              3 &
              \multicolumn{1}{c|}{3608} &
              \multicolumn{1}{c|}{4071} &
              32 &
              \multicolumn{1}{c|}{3946} &
              \multicolumn{1}{c|}{4268} &
              23 &
              \multicolumn{1}{c|}{3683} &
              \multicolumn{1}{c|}{\textbf{3904}} &
              35 \\ \cline{2-12} 
             &
              L &
              4 &
              \multicolumn{1}{c|}{3388} &
              \multicolumn{1}{c|}{3689} &
              70 &
              \multicolumn{1}{c|}{3567} &
              \multicolumn{1}{c|}{3845} &
              37 &
              \multicolumn{1}{c|}{\textbf{3352}} &
              \multicolumn{1}{c|}{\textbf{3676}} &
              42 \\ \cline{2-12} 
             &
              L &
              5 &
              \multicolumn{1}{c|}{3532} &
              \multicolumn{1}{c|}{3894} &
              45 &
              \multicolumn{1}{c|}{3582} &
              \multicolumn{1}{c|}{3926} &
              35 &
              \multicolumn{1}{c|}{\textbf{3403}} &
              \multicolumn{1}{c|}{\textbf{3751}} &
              49 \\ \cline{2-12} 
             &
              M &
              1 &
              \multicolumn{1}{c|}{3548} &
              \multicolumn{1}{c|}{3686} &
              61 &
              \multicolumn{1}{c|}{3604} &
              \multicolumn{1}{c|}{3849} &
              37 &
              \multicolumn{1}{c|}{\textbf{3478}} &
              \multicolumn{1}{c|}{\textbf{3685}} &
              66 \\ \cline{2-12} 
             &
              M &
              2 &
              \multicolumn{1}{c|}{3382} &
              \multicolumn{1}{c|}{3690} &
              52 &
              \multicolumn{1}{c|}{3203} &
              \multicolumn{1}{c|}{3798} &
              40 &
              \multicolumn{1}{c|}{3350} &
              \multicolumn{1}{c|}{\textbf{3652}} &
              37 \\ \cline{2-12} 
             &
              M &
              3 &
              \multicolumn{1}{c|}{3642} &
              \multicolumn{1}{c|}{3792} &
              55 &
              \multicolumn{1}{c|}{3490} &
              \multicolumn{1}{c|}{3821} &
              47 &
              \multicolumn{1}{c|}{3498} &
              \multicolumn{1}{c|}{\textbf{3739}} &
              68 \\ \cline{2-12} 
             &
              M &
              4 &
              \multicolumn{1}{c|}{3516} &
              \multicolumn{1}{c|}{3756} &
              45 &
              \multicolumn{1}{c|}{3588} &
              \multicolumn{1}{c|}{3760} &
              39 &
              \multicolumn{1}{c|}{\textbf{3467}} &
              \multicolumn{1}{c|}{\textbf{3651}} &
              54 \\ \cline{2-12} 
             &
              M &
              5 &
              \multicolumn{1}{c|}{3718} &
              \multicolumn{1}{c|}{3969} &
              46 &
              \multicolumn{1}{c|}{3573} &
              \multicolumn{1}{c|}{4037} &
              40 &
              \multicolumn{1}{c|}{3584} &
              \multicolumn{1}{c|}{\textbf{3804}} &
              57 \\ \cline{2-12} 
             &
              S &
              1 &
              \multicolumn{1}{c|}{4064} &
              \multicolumn{1}{c|}{4328} &
              39 &
              \multicolumn{1}{c|}{4049} &
              \multicolumn{1}{c|}{4222} &
              42 &
              \multicolumn{1}{c|}{\textbf{3917}} &
              \multicolumn{1}{c|}{\textbf{4194}} &
              37 \\ \cline{2-12} 
             &
              S &
              2 &
              \multicolumn{1}{c|}{4025} &
              \multicolumn{1}{c|}{4331} &
              41 &
              \multicolumn{1}{c|}{4042} &
              \multicolumn{1}{c|}{4459} &
              35 &
              \multicolumn{1}{c|}{\textbf{3883}} &
              \multicolumn{1}{c|}{\textbf{4202}} &
              43 \\ \cline{2-12} 
             &
              S &
              3 &
              \multicolumn{1}{c|}{3983} &
              \multicolumn{1}{c|}{4210} &
              56 &
              \multicolumn{1}{c|}{4259} &
              \multicolumn{1}{c|}{4573} &
              27 &
              \multicolumn{1}{c|}{4004} &
              \multicolumn{1}{c|}{4287} &
              37 \\ \cline{2-12} 
             &
              S &
              4 &
              \multicolumn{1}{c|}{4129} &
              \multicolumn{1}{c|}{4342} &
              38 &
              \multicolumn{1}{c|}{4091} &
              \multicolumn{1}{c|}{4546} &
              24 &
              \multicolumn{1}{c|}{\textbf{4004}} &
              \multicolumn{1}{c|}{\textbf{4254}} &
              39 \\ \cline{2-12} 
             &
              S &
              5 &
              \multicolumn{1}{c|}{3930} &
              \multicolumn{1}{c|}{4223} &
              45 &
              \multicolumn{1}{c|}{3828} &
              \multicolumn{1}{c|}{4343} &
              36 &
              \multicolumn{1}{c|}{\textbf{3778}} &
              \multicolumn{1}{c|}{\textbf{4044}} &
              48 \\ \hline
            \multirow{15}{*}{100} &
              L &
              1 &
              \multicolumn{1}{c|}{4596} &
              \multicolumn{1}{c|}{4787} &
              152 &
              \multicolumn{1}{c|}{4730} &
              \multicolumn{1}{c|}{4942} &
              158 &
              \multicolumn{1}{c|}{4667} &
              \multicolumn{1}{c|}{4877} &
              175 \\ \cline{2-12} 
             &
              L &
              2 &
              \multicolumn{1}{c|}{4283} &
              \multicolumn{1}{c|}{4875} &
              87 &
              \multicolumn{1}{c|}{4348} &
              \multicolumn{1}{c|}{4757} &
              148 &
              \multicolumn{1}{c|}{\textbf{4232}} &
              \multicolumn{1}{c|}{\textbf{4487}} &
              215 \\ \cline{2-12} 
             &
              L &
              3 &
              \multicolumn{1}{c|}{4614} &
              \multicolumn{1}{c|}{4863} &
              188 &
              \multicolumn{1}{c|}{4500} &
              \multicolumn{1}{c|}{4773} &
              227 &
              \multicolumn{1}{c|}{4535} &
              \multicolumn{1}{c|}{4904} &
              207 \\ \cline{2-12} 
             &
              L &
              4 &
              \multicolumn{1}{c|}{4195} &
              \multicolumn{1}{c|}{4698} &
              119 &
              \multicolumn{1}{c|}{4357} &
              \multicolumn{1}{c|}{4682} &
              192 &
              \multicolumn{1}{c|}{\textbf{4145}} &
              \multicolumn{1}{c|}{\textbf{4570}} &
              165 \\ \cline{2-12} 
             &
              L &
              5 &
              \multicolumn{1}{c|}{4560} &
              \multicolumn{1}{c|}{4822} &
              145 &
              \multicolumn{1}{c|}{4677} &
              \multicolumn{1}{c|}{4905} &
              187 &
              \multicolumn{1}{c|}{4702} &
              \multicolumn{1}{c|}{5173} &
              95 \\ \cline{2-12} 
             &
              M &
              1 &
              \multicolumn{1}{c|}{4738} &
              \multicolumn{1}{c|}{5054} &
              117 &
              \multicolumn{1}{c|}{4773} &
              \multicolumn{1}{c|}{5100} &
              154 &
              \multicolumn{1}{c|}{\textbf{4629}} &
              \multicolumn{1}{c|}{\textbf{4987}} &
              198 \\ \cline{2-12} 
             &
              M &
              2 &
              \multicolumn{1}{c|}{4835} &
              \multicolumn{1}{c|}{5140} &
              129 &
              \multicolumn{1}{c|}{4756} &
              \multicolumn{1}{c|}{5079} &
              201 &
              \multicolumn{1}{c|}{\textbf{4748}} &
              \multicolumn{1}{c|}{\textbf{4989}} &
              222 \\ \cline{2-12} 
             &
              M &
              3 &
              \multicolumn{1}{c|}{5058} &
              \multicolumn{1}{c|}{5374} &
              96 &
              \multicolumn{1}{c|}{4981} &
              \multicolumn{1}{c|}{5240} &
              178 &
              \multicolumn{1}{c|}{\textbf{4898}} &
              \multicolumn{1}{c|}{5304} &
              166 \\ \cline{2-12} 
             &
              M &
              4 &
              \multicolumn{1}{c|}{4827} &
              \multicolumn{1}{c|}{5174} &
              102 &
              \multicolumn{1}{c|}{4847} &
              \multicolumn{1}{c|}{5127} &
              138 &
              \multicolumn{1}{c|}{\textbf{4772}} &
              \multicolumn{1}{c|}{\textbf{5055}} &
              171 \\ \cline{2-12} 
             &
              M &
              5 &
              \multicolumn{1}{c|}{4711} &
              \multicolumn{1}{c|}{5152} &
              111 &
              \multicolumn{1}{c|}{4850} &
              \multicolumn{1}{c|}{5183} &
              150 &
              \multicolumn{1}{c|}{4869} &
              \multicolumn{1}{c|}{5319} &
              140 \\ \cline{2-12} 
             &
              S &
              1 &
              \multicolumn{1}{c|}{4735} &
              \multicolumn{1}{c|}{5000} &
              149 &
              \multicolumn{1}{c|}{4874} &
              \multicolumn{1}{c|}{5086} &
              112 &
              \multicolumn{1}{c|}{\textbf{4682}} &
              \multicolumn{1}{c|}{5111} &
              132 \\ \cline{2-12} 
             &
              S &
              2 &
              \multicolumn{1}{c|}{5139} &
              \multicolumn{1}{c|}{5552} &
              129 &
              \multicolumn{1}{c|}{5230} &
              \multicolumn{1}{c|}{5488} &
              130 &
              \multicolumn{1}{c|}{\textbf{5117}} &
              \multicolumn{1}{c|}{5503} &
              128 \\ \cline{2-12} 
             &
              S &
              3 &
              \multicolumn{1}{c|}{4867} &
              \multicolumn{1}{c|}{5134} &
              126 &
              \multicolumn{1}{c|}{5037} &
              \multicolumn{1}{c|}{5325} &
              107 &
              \multicolumn{1}{c|}{4894} &
              \multicolumn{1}{c|}{5356} &
              93 \\ \cline{2-12} 
             &
              S &
              4 &
              \multicolumn{1}{c|}{4936} &
              \multicolumn{1}{c|}{5170} &
              157 &
              \multicolumn{1}{c|}{4969} &
              \multicolumn{1}{c|}{5221} &
              152 &
              \multicolumn{1}{c|}{5013} &
              \multicolumn{1}{c|}{5310} &
              135 \\ \cline{2-12} 
             &
              S &
              5 &
              \multicolumn{1}{c|}{5106} &
              \multicolumn{1}{c|}{5294} &
              130 &
              \multicolumn{1}{c|}{4894} &
              \multicolumn{1}{c|}{5295} &
              209 &
              \multicolumn{1}{c|}{5070} &
              \multicolumn{1}{c|}{5316} &
              167 \\ \hline
            \end{tabular}
            }
    \end{table}

In columns `Best UB' and `Avg UB', the numbers in bold indicate that the corresponding value is better than the respective result of the compared methods. As observed, our refined LNS algorithm improves the best-found solution for 18 out of 30 instances - 10 for 75 requests and 8 for 100 requests. Regarding the average objective value of the best-found solutions per restart, rLNS provides the best-of-three value for 18 out of 30 instances - 13 for 75 requests and 5 for 100 requests.

The experimental results highlight the effectiveness and robustness of the proposed refined LNS (rLNS) algorithm for solving large-scale PDPT instances. Its consistent performance across diverse scenarios suggests strong generalisation capabilities, despite its reduced reliance on extensive parameter tuning. 
This makes rLNS easier to apply and more likely to perform well on unseen problem instances.

\section{Concluding remarks} \label{sec:conclusions}
This paper proposes exact and metaheuristic approaches for the Pickup and Delivery Problem with Transfers and Time Windows (PDPT), a highly applicable yet underexplored variant that generalises all classical routing problems. The exact method is driven by a novel partitioning strategy that constructs parallel paths for all pickup and delivery requests, to be then consolidated into a unified routing plan. This idea enables the design of a Logic-Based Benders Decomposition (LBBD) scheme. To address the limited scalability of the exact method, we also propose a Large Neighborhood Search (LNS) algorithm that builds on foundational ideas from related work, resulting in a more generalisable and effective solution approach.
    
The LBBD scheme, implemented as a Branch-and-Check algorithm, demonstrates significant improvements over state-of-the-art exact methods in the literature for the same problem. Experimental results on existing benchmarks show either faster convergence to optimality or reduced optimality gaps for instances with up to 30 requests, i.e., the largest available in the current literature. The performance of the LBBD approach is further enhanced when the solution from the proposed LNS algorithm is used as a warm-start.

Due to the lack of larger benchmarks, we also propose a random instance generator that incorporates more realistic assumptions. Existing works often assume large-sized requests that severely limits the number of requests accommodated in a vehicle. In contrast, our generator allows for more practical scenarios where up to 15 requests fit in a single vehicle. For these generated instances, we use the exact method to compute lower bounds, enabling a meaningful evaluation of the solutions produced by our LNS. Subsequently, we compare our LNS against state-of-the-art metaheuristic algorithms, thus showing that our method yields improved best-found solutions in most tested instances.

Our findings could motivate further effort in multiple directions. One may focus on 
more sophisticated conflict analysis, performed at solutions that lead to infeasible subproblems, to identify infeasible sub-routes; this might generate stronger feasibility cuts, further enhancing the efficiency of the method. Additionally, a hybrid matheuristic approach could be explored, where the master problem serves as a neighborhood constructor, while the LNS algorithm acts as a fast, adaptive mechanism to explore the defined neighborhood and identify improved solutions.
    
A key limitation in scaling current metaheuristic approaches for the PDPT lies in the computational cost of the `request insertion' operation, which increases significantly with the number of requests. One potential way to address this challenge is through decomposition techniques - such as those explored in \cite{Santini23} for classical vehicle routing problems - which divide the original problem into smaller, more manageable subproblems. However, applying such techniques to the PDPT is not straightforward due to the interdependence between vehicles introduced by the use of transshipment points. Designing a decomposition strategy that accounts for these interactions presents a non-trivial challenge and offers an interesting direction for future research aimed at improving scalability. 

Several routing problems are reducible to PDPT and hence eligible for applying our approach. An indicative example is the VRP with Transfers having only pickups or only deliveries. Specific problems involving mobile depots that facilitate load transfers also fall within the scope of applicability. Beyond routing, synchronisation constraints similar to those in the PDPT frequently arise in scheduling problems, hence the decomposition strategies and metaheuristic techniques presented in this paper may offer valuable insights.\\

\noindent \textbf{Acknowledgments} This research has been supported by GREEN-LOG Horizon Europe Project [grant agreement 101069892].
\setstretch{1.1}

\pagebreak
\section*{Appendix}

 \scriptsize
\begin{longtable}[c]{|c|c|c|c|cccc|cccc|}
        \caption{Exact methods on instances of \cite{Lyu23} - loose time windows}
        \label{tab:exactAll}  \\
        \hline
        \multirow{2}{*}{$|\mathcal{R}|$} &
          \multirow{2}{*}{$|\mathcal{K}|$} &
          \multirow{2}{*}{$|\mathcal{T}|$} &
          \multirow{2}{*}{Variant} &
          \multicolumn{4}{c|}{MILP \cite{Lyu23}} &
          \multicolumn{4}{c|}{LBBD} \\ \cline{5-12} 
         &
           &
           &
           &
          \multicolumn{1}{c|}{LB} &
          \multicolumn{1}{c|}{UB} &
          \multicolumn{1}{c|}{Gap (\%)} &
          Time (s) &
          \multicolumn{1}{c|}{LB} &
          \multicolumn{1}{c|}{UB} &
          \multicolumn{1}{c|}{Gap (\%)} &
          Time (s) \\ \hline
        \endfirsthead
        \endhead
        \multirow{30}{*}{10} &
          \multirow{20}{*}{2} &
          \multirow{10}{*}{1} &
          0 &
          \multicolumn{1}{c|}{649.36} &
          \multicolumn{1}{c|}{649.36} &
          \multicolumn{1}{c|}{0.00} &
          2 &
          \multicolumn{1}{c|}{649.36} &
          \multicolumn{1}{c|}{649.36} &
          \multicolumn{1}{c|}{0.00} &
          5 \\ \cline{4-12} 
         &
           &
           &
          1 &
          \multicolumn{1}{c|}{716.40} &
          \multicolumn{1}{c|}{716.40} &
          \multicolumn{1}{c|}{0.00} &
          3 &
          \multicolumn{1}{c|}{716.40} &
          \multicolumn{1}{c|}{716.40} &
          \multicolumn{1}{c|}{0.00} &
          2 \\ \cline{4-12} 
         &
           &
           &
          2 &
          \multicolumn{1}{c|}{660.26} &
          \multicolumn{1}{c|}{660.26} &
          \multicolumn{1}{c|}{0.00} &
          22 &
          \multicolumn{1}{c|}{660.26} &
          \multicolumn{1}{c|}{660.26} &
          \multicolumn{1}{c|}{0.00} &
          3 \\ \cline{4-12} 
         &
           &
           &
          3 &
          \multicolumn{1}{c|}{673.13} &
          \multicolumn{1}{c|}{673.13} &
          \multicolumn{1}{c|}{0.00} &
          12 &
          \multicolumn{1}{c|}{673.13} &
          \multicolumn{1}{c|}{673.13} &
          \multicolumn{1}{c|}{0.00} &
          6 \\ \cline{4-12} 
         &
           &
           &
          4 &
          \multicolumn{1}{c|}{821.16} &
          \multicolumn{1}{c|}{821.16} &
          \multicolumn{1}{c|}{0.00} &
          5 &
          \multicolumn{1}{c|}{821.16} &
          \multicolumn{1}{c|}{821.16} &
          \multicolumn{1}{c|}{0.00} &
          4 \\
\hline
    \end{longtable}
    \tiny
    \begin{longtable}[c]{|c|c|c|c|cccc|cccc|cccc|}
\caption{LBBD and LNS on instances of \cite{Lyu23} - loose time windows}
\label{tab:detLBBDws}\\
\hline
\multirow{2}{*}{$|\mathcal{R}|$} &
  \multirow{2}{*}{$|\mathcal{K}|$} &
  \multirow{2}{*}{$|\mathcal{T}|$} &
  \multirow{2}{*}{Variant} &
  \multicolumn{4}{c|}{LBBD} &
  \multicolumn{4}{c|}{LNS} &
  \multicolumn{4}{c|}{LBBD$^+$} \\ \cline{5-16} 
 &
   &
   &
   &
  \multicolumn{1}{c|}{LB} &
  \multicolumn{1}{c|}{UB} &
  \multicolumn{1}{c|}{Gap} &
  Time &
  \multicolumn{1}{c|}{UB} &
  \multicolumn{1}{c|}{Err} &
  \multicolumn{1}{c|}{Gap} &
  Time &
  \multicolumn{1}{c|}{LB} &
  \multicolumn{1}{c|}{UB} &
  \multicolumn{1}{c|}{Gap} &
  Time \\ \hline
\endfirsthead
\endhead
\multirow{30}{*}{20} &
  \multirow{20}{*}{2} &
  \multirow{10}{*}{1} &
  0 &
  \multicolumn{1}{c|}{973.38} &
  \multicolumn{1}{c|}{973.38} &
  \multicolumn{1}{c|}{0.00} &
  392 &
  \multicolumn{1}{c|}{999.38} &
  \multicolumn{1}{c|}{2.60} &
  \multicolumn{1}{c|}{2.60} &
  1 &
  \multicolumn{1}{c|}{973.38} &
  \multicolumn{1}{c|}{973.38} &
  \multicolumn{1}{c|}{0.00} &
  228 \\ \cline{4-16} 
 &
   &
   &
  1 &
  \multicolumn{1}{c|}{1295.13} &
  \multicolumn{1}{c|}{1295.13} &
  \multicolumn{1}{c|}{0.00} &
  189 &
  \multicolumn{1}{c|}{1409.12} &
  \multicolumn{1}{c|}{8.09} &
  \multicolumn{1}{c|}{8.09} &
  0 &
  \multicolumn{1}{c|}{1295.13} &
  \multicolumn{1}{c|}{1295.13} &
  \multicolumn{1}{c|}{0.00} &
  127 \\ \cline{4-16} 
 &
   &
   &
  2 &
  \multicolumn{1}{c|}{1012.07} &
  \multicolumn{1}{c|}{1012.07} &
  \multicolumn{1}{c|}{0.00} &
  159 &
  \multicolumn{1}{c|}{1048.37} &
  \multicolumn{1}{c|}{3.46} &
  \multicolumn{1}{c|}{3.46} &
  1 &
  \multicolumn{1}{c|}{1012.07} &
  \multicolumn{1}{c|}{1012.07} &
  \multicolumn{1}{c|}{0.00} &
  84 \\ \cline{4-16} 
 &
   &
   &
  3 &
  \multicolumn{1}{c|}{945.96} &
  \multicolumn{1}{c|}{945.96} &
  \multicolumn{1}{c|}{0.00} &
  369 &
  \multicolumn{1}{c|}{993.72} &
  \multicolumn{1}{c|}{4.81} &
  \multicolumn{1}{c|}{4.81} &
  1 &
  \multicolumn{1}{c|}{945.96} &
  \multicolumn{1}{c|}{945.96} &
  \multicolumn{1}{c|}{0.00} &
  312 \\ \cline{4-16} 
 &
   &
   &
  4 &
  \multicolumn{1}{c|}{1403.14} &
  \multicolumn{1}{c|}{1403.14} &
  \multicolumn{1}{c|}{0.00} &
  40 &
  \multicolumn{1}{c|}{1422.30} &
  \multicolumn{1}{c|}{1.35} &
  \multicolumn{1}{c|}{1.35} &
  1 &
  \multicolumn{1}{c|}{1403.14} &
  \multicolumn{1}{c|}{1403.14} &
  \multicolumn{1}{c|}{0.00} &
  32 \\
\hline
\end{longtable}
\scriptsize
    \begin{longtable}[c]{|c|c|c|c|cccc|cccc|}
\caption{Exact methods on instances of \cite{Lyu23} - tight time windows}
\label{tab:twAll}  \\
\hline
\multirow{2}{*}{$|\mathcal{R}|$} &
  \multirow{2}{*}{$|\mathcal{T}|$} &
  \multirow{2}{*}{TW} &
  \multirow{2}{*}{Variant} &
  \multicolumn{4}{c|}{MILP \cite{Lyu23}} &
  \multicolumn{4}{c|}{LBBD} \\ \cline{5-12} 
 &
   &
   &
   &
  \multicolumn{1}{c|}{LB} &
  \multicolumn{1}{c|}{UB} &
  \multicolumn{1}{c|}{Gap} &
  Time &
  \multicolumn{1}{c|}{LB} &
  \multicolumn{1}{c|}{UB} &
  \multicolumn{1}{c|}{Gap} &
  Time \\ \hline
\endfirsthead
\endhead
\multirow{180}{*}{3} &
  4 &
  180L &
  0 &
  \multicolumn{1}{c|}{434.51} &
  \multicolumn{1}{c|}{434.51} &
  \multicolumn{1}{c|}{0.00} &
  8 &
  \multicolumn{1}{c|}{434.51} &
  \multicolumn{1}{c|}{434.51} &
  \multicolumn{1}{c|}{0.00} &
  33 \\ \cline{2-12} 
 &
  4 &
  180L &
  1 &
  \multicolumn{1}{c|}{349.91} &
  \multicolumn{1}{c|}{349.91} &
  \multicolumn{1}{c|}{0.00} &
  5 &
  \multicolumn{1}{c|}{349.91} &
  \multicolumn{1}{c|}{349.91} &
  \multicolumn{1}{c|}{0.00} &
  19 \\ \cline{2-12} 
 &
  4 &
  180L &
  2 &
  \multicolumn{1}{c|}{470.34} &
  \multicolumn{1}{c|}{470.34} &
  \multicolumn{1}{c|}{0.00} &
  12 &
  \multicolumn{1}{c|}{470.34} &
  \multicolumn{1}{c|}{470.34} &
  \multicolumn{1}{c|}{0.00} &
  26 \\ \cline{2-12} 
 &
  4 &
  180L &
  3 &
  \multicolumn{1}{c|}{371.35} &
  \multicolumn{1}{c|}{371.35} &
  \multicolumn{1}{c|}{0.00} &
  8 &
  \multicolumn{1}{c|}{371.35} &
  \multicolumn{1}{c|}{371.35} &
  \multicolumn{1}{c|}{0.00} &
  13 \\ \cline{2-12} 
 &
  4 &
  180L &
  4 &
  \multicolumn{1}{c|}{368.94} &
  \multicolumn{1}{c|}{368.94} &
  \multicolumn{1}{c|}{0.00} &
  8 &
  \multicolumn{1}{c|}{368.94} &
  \multicolumn{1}{c|}{368.94} &
  \multicolumn{1}{c|}{0.00} &
  45 \\ \cline{2-12} 
 &
  4 &
  180L &
  5 &
  \multicolumn{1}{c|}{363.06} &
  \multicolumn{1}{c|}{363.06} &
  \multicolumn{1}{c|}{0.00} &
  6 &
  \multicolumn{1}{c|}{363.06} &
  \multicolumn{1}{c|}{363.06} &
  \multicolumn{1}{c|}{0.00} &
  15 \\ \cline{2-12} 
 &
  4 &
  180L &
  6 &
  \multicolumn{1}{c|}{437.78} &
  \multicolumn{1}{c|}{437.78} &
  \multicolumn{1}{c|}{0.00} &
  2 &
  \multicolumn{1}{c|}{437.78} &
  \multicolumn{1}{c|}{437.78} &
  \multicolumn{1}{c|}{0.00} &
  1 \\ \cline{2-12} 
 &
  4 &
  180L &
  7 &
  \multicolumn{1}{c|}{435.27} &
  \multicolumn{1}{c|}{435.27} &
  \multicolumn{1}{c|}{0.00} &
  2 &
  \multicolumn{1}{c|}{435.27} &
  \multicolumn{1}{c|}{435.27} &
  \multicolumn{1}{c|}{0.00} &
  4 \\
\hline
\end{longtable}

\end{document}